\documentclass[a4paper,twoside,10pt]{amsart}

\usepackage{amssymb}
\usepackage{amsmath}
\usepackage{epsfig}

\usepackage[latin1]{inputenc} 
\usepackage[T1]{fontenc}

\newtheorem{The}{Theorem}[section]
\newtheorem{Lem}[The]{Lemma}
\newtheorem{Cor}[The]{Corollary}
\newtheorem{Pro}[The]{Proposition}

\theoremstyle{definition}
\newtheorem{definition}[The]{Definition}
\newtheorem{Exam}[The]{Example}
\newtheorem{Not}[The]{Notation}

 \theoremstyle{remark}
\newtheorem{remark}[The]{Remark}
\newtheorem{remarks}[The]{Remarks}

\def\a{\alpha}
\def\A{\mathbf A}

\def\C{\mathbf C}

\def\d{\delta}

\def\f{\phi}

\def\g{\gamma}

\def\int{\mbox{\rm int}}

\def\J{{\mathcal{J}}}
\def\l{\lambda}
\def\L{\mathbf{L}}

\def\ord{\mbox{\rm ord}}

\def\O{\mbox{\rm orb}}
\def\p{\pi}

\def\Q{\mathbf Q}
\def\R{\mathbf R}
\def\r{\rho}

\def\s{\sigma}
\def\t{\tau}

\def\y{\wedge}

\def\Z{\mathbf Z}
\def\z{\zeta}

\def\limproj{\mathop{\oalign{lim\cr\hidewidth$\longleftarrow$\hidewidth\cr}}}

\title{Motivic Poincar{\'e} series, toric singularities and logarithmic jacobian ideals}

\author{H. Cobo Pablos and  P.D. Gonz{\'a}lez P{\'e}rez}

\address{Instituto de Ciencias Matem\'aticas-CSIC-UAM-UC3M-UCM. Depto. Algebra. Facultad de Ciencias Matem\'aticas. Universidad Complutense de Madrid.
Plaza de las Ciencias 3. 28040. Madrid. Spain.}

\email{hcobopab@mat.ucm.es, pgonzalez@mat.ucm.es}

\thanks{Gonz\'alez P\'erez is supported by {\em Programa Ram\'on y Cajal}  of {\em Ministerio de Educaci\'on y Ciencia}(MEC), Spain.
Cobo Pablos is supported by a grant of Fundaci\'on Caja Madrid.
Both authors are supported by  MTM2004-08080-C02-01 grant of MEC}



\begin{document}

\maketitle

{\em Abstract}. The {\em geometric motivic Poincar\'e series} of a variety,
which was introduced by Denef and Loeser, takes into account the
classes in the Grothendieck ring of the sequence of jets of arcs in the
variety. Denef and Loeser proved that this series has a rational
form. We describe it  in the case of an affine toric variety of arbitrary dimension.
The result, which provides an explicit set of candidate poles,   is expressed in terms of  the sequence of Newton polyhedra of certain monomial ideals,
which we call {\em logarithmic jacobian ideals},
associated to the modules of
differential forms with logarithmic poles outside the torus of the toric variety.

\section*{Introduction}

Let $S$ denote an irreducible and reduced algebraic variety of
dimension $d$ defined over the field $\C$ of complex numbers. The set $H(S)$ of
formal arcs of the form, $ \mbox{Spec } \C[[t]] \rightarrow S$ can
be given the structure of scheme over $\C$ (not necessarily of
finite type). If $0 \in S$ we denote by $H(S)_0 := j_0^{-1} (0) $ the subscheme of the arc space
consisting on arcs in $H(S)$ with origin at $0$.  The set $H_k(S)$ of $k$-jets of $S$,  of the form $
\mbox{Spec } \C[t] / (t^{k+1}) \rightarrow S$, has the structure
of algebraic variety over $\C$.
By a theorem of Greenberg, the
image of the space of arcs $H(S)$ by the natural morphism of
schemes $j_k : H(S) \rightarrow H_k (S)$ which maps any arc to
its $k$-jet,
  is a constructible subset of  $ H_k (S)$.
  Notice that $j_k (H(S)) = H_k (S)$ if $S$ is smooth
but  $j_k (H(S)) \ne H_k (S)$ in general.
 Since a constructible set $W$  has an image $[W]$
in the {\em Grothendieck ring  of varieties} $K_0 (\mbox{\rm
  Var}_\C) $
it is natural to measure the singularities of $S$ by considering
the formal power series:
\begin{equation}     \label{P}
P^S_{\mathrm {geom}}  (T) := \sum_{s \geq 0} [j_s(H(S))] T^s \in
K_0 (\mbox{\rm
  Var}_\C) [[T]],
\end{equation}
which is called the {\em geometric motivic Poincar\'e series} of $S$.
Similarly, the  {\em local geometric motivic Poincar\'e series} of the germ
$(S,0)$, denoted by  $P^{(S,0)}_{\mathrm {geom}}  (T)$, is defined by replacing
 $H(S)$ by $H(S)_0$ in
the right hand side of (\ref{P}).
Denef and Loeser introduced these series, inspired by the Poincar\'e series of Serre-Oesterl\'e in arithmetic geometry (see
\cite{DL-R}). They
 proved that  the image of these series  in
the ring $ K_0 (\mbox{\rm
  Var}_\C) [\L^{-1}][[T]]$ (where $\L =[\A^1_\C]$ denotes the class of the affine line)
has a rational form  (see
\cite{DL1}).

             If  $(S, 0)$ is an analytically irreducible germ of
plane curve the series  $P^{(S,0)}_{\mathrm {geom}}  (T)$ is determined  by  the
multiplicity of $(S,0)$ (see \cite{DL2}).   For a general singular variety $S$,
the invariants of $S$ encoded by the series
$P^{(S,0)}_{\mathrm {geom}} (T)$, in particular by the denominator of its
rational form, are not well understood.
In comparison with other motivic series, as the motivic zeta functions
of a polynomial or ideal, there is not a general formula for
$P^{(S,0)}_{\mathrm {geom}} (T)$ in terms of a resolution of
singularities of $S$ (see \cite{DL-bcn}). Some positive results in
this direction have been obtained by Nicaise for a class of singularities defined in terms of the
existence of an
embedded resolution
of special type; the simplest example in this class are those hypersurface
singularities with embedded resolution obtained by one blowing up (see \cite{Nicaise2}).
Lejeune-Jalabert and Reguera \cite{LR}  have
given a formula for the local geometric motivic Poincar\'e series of a germ $(S,0)$  of
{\em normal} toric surface at its distinguished point in terms of
the Hirzebruch-Jung continued fraction describing the resolution
of singularities of the germ $(S,0)$.    
The alternative computation of the geometric motivic Poincar\'e 
series of a normal toric surface singularity is given in \cite{Nicaise2}. 
The comparison with the arithmetic and the Igusa series is contained in  \cite{Nicaise}.

In this paper we consider the case of $(S,0)$ being  the germ of an affine toric variety 
$(Z^\Lambda,0)$ of dimension $d$
at its $0$-dimensional orbit.
Here  $\Lambda$ denotes a semigroup of finite type of a rank $d$ lattice $M$ such that $Z^\Lambda :=
\mbox{{\rm Spec}} \C[\Lambda]$  (cf.
~ Notation \ref{Lambda}).

Our approach is inspired by Lejeune-Jalabert and Reguera \cite{LR}, though there are substantial differences coming from the
particularities of normal toric surfaces, which verify that:
\begin{itemize}
 \item[(a)] Every truncated arc is the jet of an arc with generic point in the torus.

\item[(b)] Any pair of consecutive integral vectors in the boundary of the Newton polygon of the maximal ideal define a basis of the lattice.
\end{itemize}
Property (a) holds more generally for normal toric singularities (see \cite{Nicaise}) but 
 does not hold in general without the assumption of normality, the simplest example is the Whitney umbrella 
(see Remark \ref{ejemplo}). Property (b), which plays  also an essential role in the comparison of various types of motivic series in \cite{Nicaise},  does not generalize even for normal toric germs of dimension $\geq 3$.

We deal with the failure of property (a) by characterizing combinatorially the jets of those 
arcs which cannot be obtained as jets of arcs factoring through proper orbit closures of the 
action of the torus on $Z^\Lambda$. We define the auxiliary  series
$P( \Lambda) $ by taking classes in the Grothendieck ring of these sets and considering the associated Poincar\'e series.
We have that
$ P_{\mathrm {geom}}^{(Z^\Lambda,0)}(T) = \sum P(\Lambda \cap \t )$,
where $\t$ runs through the faces of the cone $\R _{\geq 0} \Lambda$ of $M_\R :=  M \otimes_\Z \R$. The term
$P(\Lambda \cap \t)$ is the auxiliary series associated to the toric variety $Z^{\Lambda \cap \t}$, which is an orbit closure of the torus action
on $Z^\Lambda$.

The failure of (b) is overcome by the systematic use of the {\em logarithmic jacobian ideals} associated 
to the toric variety $Z^{\Lambda}$ to study jet spaces.  
The logarithmic jacobian ideals $\J_1, \dots, \J_d$ of $Z^\Lambda$ are defined in terms of the minimal set
 of generators of the semigroup $\Lambda$ in Section \ref{main}. The ideal $\J_1$ is the maximal ideal defining 
the closed point of the germ $(Z^{\Lambda},0)$. The ideal $\J_d$ appears in \cite{LR} in connection with 
the combinatorial description of the Nash blowing up (see also
\cite{GN, Teissier}). If $1\leq k \leq  d$, the logarithmic jacobian ideal $\J_k$ can be described 
in terms of the module of  K\"{a}hler differential $k$-forms on
$Z^\Lambda$ over $\C$, in a way which generalizes the one given for $\J_d$ in the Appendix of \cite{LR} 
(see Section \ref{tor}).
Up to our knowledge, if $d \geq 3$ the ideals $\J_2, \dots, \J_{d-1}$ appear in this paper for the first time in
the literature.

Ishii noticed that the arc space of the torus acts on the arc space of the toric variety $Z^\Lambda$ (see \cite{Ishii-algebra,
Ishii-crelle}).
The set $H^*$ consisting of those arcs of $H(Z^\Lambda)_0$ which have their generic point in the torus is  a union of orbits.
These orbits are in bijection with the possible orders of the arcs,
naturally identified with the elements of the dual lattice $N:= M^*$, which are in the interior of the dual cone $\s$ of $\R_{\geq 0} \Lambda$.
For $\nu \in \stackrel{\circ}{\s} \cap N$ we denote by $H^*_\nu$ the corresponding orbit in the arc space.
We show that the jets of these orbits are either disjoint or equal and we characterize the equality in combinatorial
terms.
We prove that the coefficient of $T^m$ in the auxiliary series $P(\Lambda)$ expands as the sum of classes $[j_m (H^*_\nu) ]$ in the Grothendieck ring, for $\nu$ running
through a finite subset of $\stackrel{\circ}{\s} \cap N$.
The combinatorial convexity properties of the Newton polyhedra of the logarithmic jacobian ideals
allow us to determine a simple formula for the class of $j_m (H^*_\nu)$ in the
Grothendieck ring (see Theorem \ref{key-bis}).

The main result states that the rational form of  the
geometric motivic Poincar\'e series $P_{\mathrm {geom}}^{(Z^\Lambda,0)}(T)$
is determined by the Newton polyhedra (with integral structure) of the logarithmic jacobian ideals
of the orbit closures of $Z^\Lambda$ (see Theorem \ref{PLambdaRac} and Corollary \ref{P-geom}).
In particular we describe explicitly a  finite set of candidate poles for the rational
form of $P_{\mathrm {geom}}^{(Z^\Lambda,0)}(T)$. We give a geometrical interpretation
of the candidate poles in terms of the order of vanishing
of certain sheaves of locally principal monomial ideals along the exceptional divisors of certain modifications,
 which are both defined in terms of the logarithmic jacobian ideals.
The rationality of the series  is deduced at this point from  a  purely combinatorial result:
the rationality of the generating series of the projection of the set of integral points in the
interior of a rational open cone  (see Theorem \ref{teoremaapendice}). The appearance of these projections
seems the combinatorial analogue of the quantifier elimination results  used in \cite{DL1}.

We give two applications:
\begin{itemize}
\item We deduce a formula for the global geometric motivic Poincar\'e series $P_{\mathrm {geom}}^{Z^\Lambda}(T)$
in the normal case (see Theorem \ref{P-locales}).

\item We prove a  formula  for the {\em motivic volume} of the arc space of the germ
$(Z^\Lambda, 0)$ in terms of the logarithmic jacobian ideal
$\mathcal J_d$ (see Proposition \ref{vol-motGlobal2}). We have
obtained this result without using Denef and Loeser's formula for
the motivic volume of a variety $S$ in terms of a resolution of
singularities (see \cite{DL1}).
\end{itemize}

In the normal toric surface  case,  property (b) allows an
explicit description of the series in \cite{LR}. In this case only
the terms $1- T$ and $1 -\L T$, which appear then in the
denominator of the rational form of the series in Corollary
\ref{P-geom}, are not actual poles. This property is a
particularity of the normal toric surface case. We
give an example   of toric surface  $Z^\Lambda$ such that all
        terms in the denominator of the rational form of
$P_{\mathrm {geom}}^{(Z^\Lambda,0)}(T)$ in Corollary \ref{P-geom},
are actually poles (see Section
\ref{toric-example}).

In \cite{C-GP-qo} we extend the results and approach of this paper
to the case of a germ of irreducible quasi-ordinary hypersurface
singularity of arbitrary dimension $d$ in terms of similar notions of logarithmic jacobian ideals.
Rond states some partial results on  this case in \cite{Rond}.
In general it is a challenge to analyse this motivic series in terms of some suitable notion of logarithmic jacobian ideals
associated to a partial resolution of singularities of a given singularity.  It is a natural problem to find
 a geometrical meaning for the logarithmic jacobian ideals in terms of
limits of tangent spaces.

The results of this paper hold if the base field of complex numbers $\C$  is replaced by 
an algebraically closed field of zero
characteristic. The assumption that the base field has characteristic zero is used in 
Section \ref{universal}.

The paper is organized as follows. In Section \ref{sec-tor} we set our notations on toric varieties.
In Section  \ref{intro-arcs} some results on arcs and jets
spaces are recalled. We describe the orbit decomposition of  the arc space of a
toric variety in Section \ref{ArcsAndJets}. In Section \ref{main}
we state the main results. In Section \ref{conv-newton} we give
some  combinatorial convexity properties of the Newton polyhedra
of the logarithmic jacobian ideals. Section \ref{universal} deals with the universal family of arcs in the torus. In Section \ref{tor-jet} we analyze the jets of the orbits in the arc space. The main results on the
geometric motivic Poincar\'e series are proved in Sections
\ref{GenptOrbits} and \ref{racionality}. A formula for the motivic
volume is given in Section \ref{volume}. Sections \ref{tor} and
\ref{app} can be read independently of the rest of the paper.
Section \ref{tor} is dedicated to  the definition of the sequence of
logarithmic jacobian ideals in terms of differential forms.
Section  \ref{app}  deals with  generating functions.

\section{A reminder of toric geometry} \label{sec-tor}

In this Section we introduce  the basic notions and notations  from toric
geometry (see \cite{Ewald,Oda,Fu,GKZ}
  for
proofs).
Following the convention established at the meeting
``Convex and algebraic geometry'',  Oberwolfach (2006), we do not assume the normality in the definition of
toric varieties.

If $N \cong \Z^{d}$ is a lattice we denote by $N_\R:=N\otimes\R$
(resp. $N_\Q:=N\otimes\Q$) the vector space spanned by $N$ over
the field $\R$ (resp. over $\Q$). In what follows a {\em cone} in
$N_\R$ mean a {\em rational convex polyhedral cone}: the set of
non negative linear combinations of vectors $a_1, \dots, a_r \in
N$. The cone $\t$ is {\em strictly convex} if it contains no line
through the origin,  in that case we denote by $0$ the
$0$-dimensional face of $\t$;
the cone $\t$ is {\em simplicial} if the primitive vectors of the
$1$-dimensional faces are linearly independent over $\R$. We
denote by $\stackrel{\circ}{\t}$ or by $\mbox{{\rm int}} (\t)$ the
relative interior of the cone $\t$. We denote by $\R \t$ the real
vector subspace spanned by $\t$ in $N_\R$.

 We denote by $M$ the dual
lattice. The {\em dual} cone  $\t^\vee \subset M_\R$ (resp. {\em
orthogonal} cone $\t^\bot$) of $\t$ is the set $ \{ w  \in M_\R\ |\
\langle w, u \rangle \geq 0,$  (resp. $ \langle w, u \rangle = 0$)
$ \; \forall u \in \t \}$.

A {\em fan} $\Sigma$ is a family of strictly convex
  cones  in $N_\R$
such that any face of such a cone is in the family and the
intersection of any two of them is a face of each. The relation
$\theta \leq \t$  (resp. $\theta < \t$) denotes that $\theta$ is a
face of $\t$ (resp. $\theta \ne \t$ is a face of $\t$). The {\em
support} (resp. the $k$-{\em skeleton}) of the fan $\Sigma$ is the
set $|\Sigma | := \bigcup_{\t \in \Sigma} \t \subset N_\R$ (resp.
$\Sigma^{(k)} = \{ \t \in \Sigma \mid \dim \t = k \}$).
We say that a fan $\Sigma'$ is a {\em subdivision\index{fan
subdivision}} of the fan $\Sigma$ if both fans have the same
support and if every cone of $\Sigma'$ is contained in a cone of
$\Sigma$. If $\Sigma_i$ for $i=1, \dots, n$, are fans with the
same support their intersection
$\cap_{i=1}^n\Sigma_i:=\{\cap_{i=1}^n\t_i\ |\
\t_i\in\Sigma_i\}$ is also a fan.   The $1$-skeleton of  $\cap_{i=1}^n\Sigma_i$ is    $\cup_{i=1}^n\Sigma_i^{(1)}$.

\begin{Not} \label{Lambda}
In this paper $ \Lambda  $ is a sub-semigroup of finite type of a
lattice $M$,  which generates $M$ as a group and such that the
cone $\s^\vee = \R_{\geq 0} \Lambda$ is strictly convex and of
dimension $d$. We denote by $N$ the dual lattice of $M$ and by $\s
\subset N_{\R}$ the dual cone of $\s^\vee$. We denote by
$Z^\Lambda$ the {\em affine toric variety} $Z^{\Lambda} =
\makebox{Spec} \, \C [ \Lambda ]$, where $ \C[ \Lambda] = \{
\sum_{\rm finite}  a_\l {X}^\l \mid a_{\l} \in \C \} $ denotes the
semigroup algebra of the semigroup $\Lambda$ with coefficients in
the field $\C$. The semigroup $\Lambda$ has a unique minimal set
of generators $e_1,\dots,e_n$ (see the proof of  Chapter V, Lemma
3.5, page 155 \cite{Ewald}). We have an {\em  embedding} of
$Z^\Lambda\subset \C^n$ given by, $x_i := X^{e_i}$ for $i=1, \dots
n$.
\end{Not}

If $\Lambda = \s^\vee \cap
M$   then the variety    $Z^{\Lambda}$,
which we denote also by         $Z_{\s, N}$ or by
 $Z_\s$ when the lattice is clear from the context, is normal.
If $\Lambda \ne \s^\vee \cap
M$       the
inclusion of semigroups $\Lambda \rightarrow \bar{\Lambda}:= \s^\vee \cap M $
defines a toric modification $Z^{\bar{\Lambda}} \rightarrow
Z^{{\Lambda}}$,
  which  is the {\em normalization map}.

The
 affine varieties $Z_\t$ corresponding to cones in a fan $\Sigma$
 glue up to define a {\em toric variety\index{toric variety}} $
 Z_\Sigma$.
For instance, the toric variety defined by the fan formed by the
 faces of the cone $\s$ coincides with the affine toric variety
 $Z_\s$.
The subdivision $\Sigma'$ of  a fan $\Sigma$ defines a {\em toric
 modification\index{modification}} $ \pi_{\Sigma'} : Z_{\Sigma '}
 \rightarrow   Z_\Sigma$.

The torus $ T_N:= Z^{M}$ is an open dense subset of $Z^\Lambda$,
which acts on $Z^\Lambda$ and the action extends the action of the
torus on itself by multiplication.
The {\em origin} $0$ of the affine toric variety $Z^{\Lambda}$
is the $0$-dimensional orbit, which is defined by the maximal ideal
$(X^{\l})_{0 \ne \l \in \Lambda }$ of $\C [\Lambda]$. There is a one to one inclusion
reversing correspondence between the faces  of $\s$ and the orbit
closures  of the torus action on $Z^\Lambda$.
If $\theta \leq \s$, we denote by  $\O^\Lambda_\theta$ the orbit
corresponding to the face $\theta$ of $\s$. The set $\Lambda \cap
\theta^\bot$ is a semigroup of finite type which generates a
sublattice $M(\theta,\Lambda)$ of finite index $i( \theta,
\Lambda)$ of the lattice $M \cap \theta^\bot$. We denote by
$N(\theta,\Lambda)$ the dual lattice of $M(\theta,\Lambda)$. The
image $\s/ \R \theta$ of the cone $\s$ in the vector space $N_\R /
\R \theta$ is the dual cone of $\R_{\geq 0} (\Lambda \cap
\theta^\bot)$.  The toric variety $Z^{\Lambda \cap
\theta^\bot }$, which is embedded in $Z^\Lambda$, is the closure
of $\O^\Lambda_\theta$. The origins of $Z^\Lambda$
and $Z^{\Lambda \cap \theta^\bot }$ coincide.
We say that  $Z^\Lambda$ is {\em analytically unibranched} if for
any $x \in Z^\Lambda$ the germ $(Z^\Lambda, x)$ is analytically
irreducible, i.e., for all $\theta \leq  \s$ we have $i(\theta,
\Lambda) =1$ (see \cite{GKZ} Chapter 5). Notice that normal toric
varieties are analytically unibranched but the converse is not
true.

The ring $\C [[
\Lambda ]] $ of formal power series with coefficients in $\C$ and
exponents in the semigroup $\Lambda$ is isomorphic to the
completion of the local ring of germs of holomorphic functions at
$(Z^\Lambda, 0)$ with respect to its maximal ideal.

The {\em Newton polyhedron} of a monomial ideal corresponding to a
non empty set of lattice vectors
 ${\mathcal I} \subset   \Lambda $ is defined as the convex hull
of the Minkowski sum  of sets ${\mathcal I} + \s^\vee$. We denote
this polyhedron by ${\mathcal N} ({\mathcal I})$.
We denote by $\mbox{\rm ord}_{\mathcal {\mathcal I}}$ the {\em support
function} of the polyhedron ${\mathcal N} ( {\mathcal I} )$, which
is defined by $\mbox{\rm ord}_{\mathcal I} : \s \rightarrow \R$, $
\nu \mapsto \inf_{\omega \in {\mathcal N} ( {\mathcal I} )}
\langle \nu, \omega \rangle$.
The face  of  the polyhedron ${\mathcal N} ({\mathcal I})$
determined by $\nu \in \s$ is  the set $ {\mathcal F}_\nu := \{
\omega \in {\mathcal N} (\mathcal I) \mid \langle \nu, \omega
  \rangle   = \mbox{\rm ord}_{\mathcal I} (\nu ) \}$.
All faces of ${\mathcal N} ({\mathcal I})$ are of this form, the
compact faces are defined by vectors $\nu  \in
\stackrel{\circ}{\sigma}$. The {\em dual fan}   $\Sigma ({\mathcal
I})$ associated to an integral polyhedron ${\mathcal N} (\mathcal I)$ is a
fan supported on $\s$ which is  formed by the cones $ \s(
{\mathcal F} ) := \{ \nu \in \s \; \mid \langle \nu , \omega
\rangle   = \mbox{\rm ord}_{\mathcal I} ( \nu ),
 \; \forall \omega \in {\mathcal F}\}$,
for ${\mathcal F}$ running through the faces of ${\mathcal N}
({\mathcal I})$. Notice that if $\t\in  \Sigma ( {\mathcal I} )$ and if
$\nu, \nu' \in \stackrel{\circ}{\t}$ then ${\mathcal
F}_{\nu} = { \mathcal F}_{\nu'}$. We denote this face of
${\mathcal N} (\mathcal I)$ also by  ${ \mathcal F}_{\t}$. Notice
that the vertices of ${\mathcal N} (\mathcal I)$ are elements of
$\mathcal I$.

If $\mathcal{I}$ is a monomial ideal of $Z_\s$ and  $\Sigma = \Sigma ({\mathcal I}) $,  then the toric modification
  $ \p_\Sigma : Z_{\Sigma}   \rightarrow Z_\s$
  is the {\em normalized blowing up} of $Z_\s$ centered at ${\mathcal I}$ (see \cite{LR} for instance).

\section{Arcs, jet spaces and the geometric motivic Poincar\'e series} \label{intro-arcs}

In this Section we introduce arc and jet spaces on a variety  $S$, i.e., a reduced separated scheme of finite type over 
$\C$.  For simplicity 
we assume that $S$ is  affine and equidimensional of dimension $d$. 
We refer to
\cite{Ishii-sedano, EM} for expository papers on arc
and jet schemes. See
\cite{DL-bcn, Loo, Veys}, for expository papers on motivic
integration on arc spaces and applications.

We have that for all integers $m \geq 0$ the functor from the
category of $\C$-algebras to the category of sets, sending a
$\C$-algebra $R$ to the set of $R[t]/(t^{m+1})$-rational points of
$S$ is representable by a $\C$-scheme $H_m(S)$ of finite type over
$\C$, called the $m$-jet scheme of $S$. The natural maps induced
by truncation $j_m^{m+1}: H_{m+1} (S) \rightarrow H_m (S)$ are
affine and hence the projective limit $H(S):= \limproj H_m (S)$ is
a $\C$-scheme, not necessarily of finite type, called the arc
space of $S$. The scheme $H(S)$ represents a functor sending a
$\C$-algebra $R$ to the set of $R[[t]]$-rational points of $S$. It
is the arc space of $S$. We consider the schemes $H_m(S)$ and
$H(S)$ with their reduced structure.      If $Z \subset S$ is a closed subvariety then $H(S)_Z:= j_0 ^{-1}
(Z) $ (resp. $H_{m}(S) _Z:= (j^m_0 )^{-1} (Z)$) denotes the
subscheme of $H(S)$ (resp. of $H_{m}(S)$) formed by arcs (resp. $m$-jets) in $S$
with origin in $Z$.

 We have natural morphisms
$j_m: H(S) \rightarrow H_m (S)$. By an {\em arc} we mean a
$\C$-rational point of $H(S)$, i.e., a morphism $ \mbox{Spec }
\C[[t]] \rightarrow S$. By an $m$-jet we mean a $\C$-rational
point of $H_m(S)$, i.e., a morphism $ \mbox{Spec } \C[t]/(t^{m+1})
\rightarrow S$. The origin of the arc (resp. of the $m$-jet) is
the image of the closed point $0$ of $ \mbox{Spec } \C[[t]]$
(resp. of $ \mbox{Spec } \C[t]/(t^{m+1}) $).

If $h(t)  =  \sum_{i \geq 0} a_i t^i $ is a formal power series and $m \geq 0$ we set
$j_m(h(t) ) :=  h(t) \mod t^{m+1}$.

Suppose that  $S\subset \A^n_\C$  is a closed affine subvariety  with ideal $I
\subset \C[x_1, \dots, x_n]$,  for $(x_1, \dots, x_n)$ coordinates
of $\C^n$.  An arc $\mbox{Spec } \C[[t]] \rightarrow \A^n_\C$
 is defined by $n$
power series
\begin{equation} \label{arc-exp}
x_i (t) =  a_i^{(0)} + a_i^{(1)} t + a_i^{(2)}t^2 + \cdots +
a_i^{(r)}t^r + \cdots      , \quad  i=1,\dots, n.
\end{equation}
 An $m$-jet  $\mbox{Spec} \C [t] /(t^{m+1})  \rightarrow \A^n_\C$
is defined by $n$-polynomials of the form given by (\ref{arc-exp}) $\mod t^{m+1}$.
If $F \in I$ we have a power series expansion
\begin{equation}  \label{arc-condition}
 F(x_1(t), \dots, x_n (t)) = \a_F^{(0)} ( \underline{a}^{(0)} ) + \a_F^{(1)} ( \underline{a}^{(0)} , \underline{a}^{(1)}
 )t +  \a_F^{(2)} ( \underline{a}^{(0)} , \underline{a}^{(1)}, \underline{a}^{(2)}
 )t^2 +\cdots
\end{equation}
where the coefficients $\a_F^{(k)}$ are polynomials expressions in $(
\underline{a}^{(0)} , \dots,  \underline{a}^{(k)})$ where
$\underline{a}^{(j)} = ( a^{(j)}_1, \dots, a^{(j)}_n )$, for $j
\in \Z_{\geq 0}$.
The arc (\ref{arc-exp}) (resp. the $m$-jet of (\ref{arc-exp}))
factors through $S$ if  (\ref{arc-condition}) vanishes for all $F \in I$
(resp. if (\ref{arc-condition}) vanishes  $\mod t^{m+1}$ for all $F \in I$).
The arc space $H(S)$  (resp.     $m$-jet space $H_m (S)$)
is the reduced scheme
underlying the affine scheme $\mbox{Spec}  \mathcal{A}_{S} $,  where
$\mathcal{A}_{S} = \C [
\underline{a}^{(0)} , \dots, \underline{a}^{(k)}, \dots ]
/(\a_F^{(0)}, \dots, \a_F^{(k)}, \dots)_{F \in I}$
(resp.            $\mbox{ Spec } \mathcal{A}_{S,m}$, where
$
\mathcal{A}_{S,m} := \C [ \underline{a}^{(0)} , \dots,
\underline{a}^{(m)}] / ( \a_F ^{(0)}, \dots, \a_F ^{(m)})_{F \in
I}$).
The {\em universal family of arcs} of $S$, which  is  the map $
\mbox{ Spec }  \mathcal{A}_{S} [[t]] \rightarrow S $ defined by
(\ref{arc-exp}), parametrizes the arcs in $H(S)$.

We recall the definition of the {\em Grothendieck ring} $ K_0
(\mbox{\rm
  Var}_\C) $
of $\C$-{\em varieties}. This ring is generated by the symbols
$[X]$ for $X$ an algebraic variety,  subject to relations: $[X] =
[X']$  if $X$ is isomorphic to $X'$, $[X] = [X- X'] + [X'] $ if
$X'$ is closed in $X$ and $[X][X'] = [X \times X']$.    We denote by $\L:= [\A_{\C}^1]$ the class of
the affine line and by $\mathcal{M}$ the localization  $ K_0
(\mbox{\rm Var}_\C) [\L^{-1}] $.

If $C$ is a constructible subset of some variety $X$, i.e. a disjoint union of finitely
many locally closed subvarieties $A_i$ of $X$ , then it is easy to see that $[C] \in K_0
(\mbox{\rm
  Var}_\C) $ is well defined as $[C] :=\sum_i [A_i ]$.

A set $A \subset H(S)$ is {\em constructible} or {\em cylindric}
if $A = j_m^{-1} (C)$,  for  some integer $m$ and some
constructible subset $C \subset H_m (S)$; the constructible set
$A$ is {\em stable} if, in addition, for all $p \geq m$ the projection $j^{p+1}_p : j_{p+1} (A) \rightarrow j_p (A)$
is a piece-wise trivial fibration with fiber $\A^d_\C$ (where $d = \dim S$).  If $A \subset H(S)$ is
constructible and $A \cap H  (\mbox{\rm Sing} S) = \emptyset$ then
$A$ is stable  (see \cite{DL1}). For a stable set it makes sense
to consider the {\em naive motivic measure}, defined as the limit
${\lim}_{m \rightarrow \infty} [j_m (A) ] \L^{-md} \in
\mathcal{M}$ (by definition of stability all the terms $ [j_m (A)
] \L^{-md}$ are equal for $m$ large enough).  Kontsevich introduced a completion  $\hat{\mathcal{M}}: =
\underleftarrow{\lim} \mathcal{M}/F^m$ of the ring $\mathcal{M}$, where   $F^{m}$  for $m \in \Z$,   is the subgroup of $\mathcal{M}$
generated by $[X] \L^{-i}$ such that $\dim X + m \leq i$ and $(F^m)$ defines a ring
filtration since $F^m F^p \subset F^{m+p} $.
\begin{The} {\rm (see \cite{DL1} Theorem 7.1)}
Let $A$ be a constructible subset of $H(S)$. Then the limit $\mu
(A) := {\lim}_{m \rightarrow \infty} [j_m (A) ] \L^{-md}$ exists in $\hat{\mathcal{M}}$.
If $A=H(S)$ this limit is nonzero. \label{ThVolMot}
\end{The}
If $A$ is a constructible subset of $H(S)$, then $\mu
(A)$
is called the {\em motivic measure} of $A$.  Notice that if $S$ is irreducible and $Z \subset S$ is a proper closed subset then $H(Z)
\subset H(S)$ is not cylindric. There exists a class of {\em
measurable} sets containing $H(Z)$ and the cylinders and a measure
$\mu$ with values on $\hat{\mathcal{M}}$, extending the
motivic measure of constructible sets. We refer to \cite{DL1, DL-M, Loo}
for the precise definition.
\begin{definition}
The motivic measure of the arc space $H(S)_Z$ for $Z$ a
closed subvariety of $S$, is called the {\em motivic volume} of
$H(S)_Z$.
\end{definition}
\begin{Pro} {\rm (see \cite{DL1, DL-M, Loo})} \label{ProVolMot}
 If $A \subset H(S)$ is a measurable set such that $A
\subset H(Z)$ for some closed subvariety $Z \subset S$ with $\dim
Z < \dim S$ then $\mu (A) =0$.
\end{Pro}

 By a Theorem of Greenberg
\cite{Greenberg}, see also \cite{EM},  $j_m (H(S))$ is
a constructible subset of $H_m (S)$, hence it has an image in the
Grothendieck ring $ K_0 (\mbox{\rm
  Var}_\C) $. The same applies for $j_m (H(S)_Z)$ if $Z\subset S$ is a closed
  subvariety.
\begin{definition} Let $S$ be a variety and $Z \subset S$  a closed subvariety.
The {\em geometric motivic Poincar\'e series} of $S$ (resp. of $(S,Z)$)
is the element of $K_0 (\mbox{\rm  Var}_\C) [[T]]$ defined by
\[
P_{\mathrm {geom}}^{S} (T) := \sum_{m \geq 0} [ j_m (H(S)) ] T^m   \quad
\mbox{(resp. }  P_{\mathrm {geom}}^{(S, Z)} (T) := \sum_{m \geq 0} [ j_m (H(S)_Z)
] T^m \, ).
\]
\end{definition}
For instance, it is easy to see that
$P_{\mathrm {geom}}^{\C^d}(T)={\L^d}{(1-\L^dT)^{-1}}$ and $P_{\mathrm {geom}}^{(\C^d,0)}(T)={(1-\L^dT)}^{-1}$.
We often call the series $P_{\mathrm {geom}}^S(T)$ (resp. $P_{\mathrm {geom}}^{(S,Z)}(T)$)
the motivic Poincar\'e series of $S$ (resp. of $(S,Z)$) for short.
Denef and Loeser proved that these series have a rational form:

\begin{The}  {\rm  (see \cite{DL1}   Theorem 1.1)}
The series $P_{\mathrm {geom}}^S (T)$ (resp.  $P_{\mathrm {geom}}^{(S, Z)} (T)$),
considered as an element of $\mathcal{M} [[T]]$ belongs to
$\mathcal{M} (T)$, more precisely there exist $Q(T) \in
\mathcal{M}[T]$, $a_i \in \Z$ and $b_i \in \Z_{\geq 1}$, for $i=1,
\dots,r$,  such that the series is of the form $Q(T) \prod_{i=1}^r
(1 - \L^{a_i} T^{b_i})^{-1}$. \label{RacDL}
\end{The}
The proof of this deep result is based on quantifier elimination for
semi-algebraic sets of power series, a substantial development of
the theory of motivic integration introduced by Kontsevich and the
existence of resolution of singularities of varieties over a field
of zero characteristic.  See \cite{DL-R} for relations with
other Poincar\'e series in arithmetic geometry.

\section{Arcs and jets on a toric singularity}

\label{ArcsAndJets}

Let $\Lambda$ be a semigroup, as in Notation \ref{Lambda}.
If $R$ is a $\C$-algebra, a $R$-rational point
of $Z^{\Lambda}$ is a homomorphism of semigroups $(\Lambda, +)
\rightarrow (R,\cdot)$, where $(R,\cdot)$ denotes the semigroup
$R$ for the multiplication. In particular, the closed points are
obtained for $R = \C$.
An arc
$h$ on the affine toric variety $Z^{\Lambda}$ is given by a semigroup
homomorphism $(\Lambda, +) \rightarrow (\C[[t]], \cdot)$. An
arc in the torus $ T_N$ is defined by a semigroup homomorphisms $\Lambda
\rightarrow \C[[t]]^*$, where $\C[[t]]^*$ denotes the group of
units of the ring $\C[[t]]$.

\begin{Not}    We denote  the set of arcs $H(Z^\Lambda)_0$
of $Z^\Lambda$ with
origin at the distinguished point $0$ of $Z^\Lambda$  simply by $ H_\Lambda$, and by
    $H_\Lambda^*$ the set consisting of those arcs of
$H_\Lambda$ with generic point in the torus $ T_N$.
\label{NotacionArcs}
\end{Not}
Notice that $h \in H_\Lambda^*$ if and only if for all $u \in
\Lambda$ the formal power series ${X}^u \circ h \in
\C [[t]]$ is non-zero.
Any arc  $h  \in H^*_\Lambda$ defines two group
 homomorphisms $\nu_h : M \rightarrow \Z  \mbox{ and } \omega_h: M
\rightarrow \C[[t]]^* \mbox{ by: } {X}^m \circ {h} = t^{\nu_h (m)}
\omega_h (m)$.  If $m \in \Lambda$ then  $\nu_h (m)
>0$ hence $\nu_h$ belongs to $\stackrel{\circ}{\s} \cap N$.
Notice that $\omega_h$ defines an arc in the torus, i.e.,
$\omega_h \in H ({ T_N}) $.

Ishii noticed that the space of arcs in the torus acts on the arc
space of a toric variety (see \cite{Ishii-algebra, Ishii-crelle}).
\begin{Lem} \label{encode} {\rm (Theorem 4.1 of \cite{Ishii-algebra}, and Lemma 5.6 of
 \cite{Ishii-crelle}).}
The map $ \stackrel{\circ}{\s} \cap N \times H({ T_N}) \rightarrow
H^*_\Lambda$  which sends a pair $(\nu,\omega)$ to the arc $ h$
defined by
$
{X}^u \circ {h} = t^{\langle \nu, u \rangle} \omega (u), \mbox{
for } u \in \Lambda ,
$
is a one to one correspondence. The sets $H^*_{\Lambda, \nu} := \{
h \in H^*_\Lambda\ |\ \nu_h = \nu \}$   for $\nu \in
\stackrel{\circ}{\s} \cap N$ are orbits for the action of $H_{
T_N} $ on $H^*_\Lambda$ and we have that $H_\Lambda^* =
\bigsqcup_{\nu \in  \stackrel{\circ}{\s} \cap N}  H^*_{\Lambda,
\nu} $.
\end{Lem}
       The
sets defining these orbits were also considered by
Lejeune-Jalabert and Reguera in the normal toric surface case
(Proposition 3.3 of \cite{LR}).

\begin{remark}
We often denote the set $H^*_\Lambda$ (resp. the orbit
$H^*_{\Lambda, \nu} $) by  $H^*$ (resp. by $H^*_\nu$) if $\Lambda$
is clear from the context.
\end{remark}

 An arc $h \in H_{\Lambda}$
has its generic point $\eta$ contained in exactly one orbit of the
torus action on $Z^\Lambda$.  If    $h(\eta) \in \O_\theta^\Lambda$, for some $\theta \leq \s$, then
$h$ factors through the orbit closure $Z^{\Lambda
\cap \theta^\bot}$ and  $h \in H_{\Lambda \cap \theta^\bot}^*$, i.e., $h$ is an arc through $(Z^{\Lambda
\cap \theta^\bot}, 0)$ with generic point in  $T_{N(\Lambda, \theta)} = \O_\theta^\Lambda$.
     We can apply Lemma \ref{encode} to describe the set
$H_{\Lambda \cap \theta^\bot}^*$, just replacing $\Lambda$, $\s$,
$M$ and $N$ by $\Lambda \cap \theta^\bot$, $\s / \R \theta  $, $M
(\Lambda, \theta)$ and $N(\Lambda, \theta)$ respectively (cf. with
notations in Section \ref{sec-tor}).
 In
particular, if $\theta = 0$ then $h \in H_\Lambda^*$;  if $\theta
=\s$ then $\Lambda \cap \theta^\bot =0$ and $h$ is the constant
arc at the distinguished point $0 \in Z^\Lambda$.
We have a
partition $     H_\Lambda = \sqcup_{\theta \leq \s}  H_{\Lambda \cap
    \theta^\bot}^*$.

\begin{remark}
In the normal case the equality $j_m(H_\Lambda)=j_m(H_\Lambda^*)$ holds for all $m \geq 0$,
see \cite{Nicaise}. This property fails in general, for
instance, the arc $h(t) = (0, t, 0)$ of the {\em Whitney umbrella}, $\{ (x_1, x_2, x_2) \mid x_1 ^2 x_2 -x_3^2 =0 \}$, is contained in
the singular locus but its $1$-jet is not obtained as the jet of
an arc $h'$ with generic point in the torus.
        \label{ejemplo}
\end{remark}

\section{Statement of the  main results on the geometric motivic Poincar\'e
series} \label{main}

In this Section we state the main results of the paper. The proofs are given in  Section \ref{racionality}.

We
consider the following
auxiliary Poincar\'e series:

\begin{equation}
   P(\Lambda):=
    \displaystyle\sum_{s\geq 0}\big[j_s(H_\Lambda^*)\setminus\displaystyle\bigcup_{0 \neq\theta\leq\s}
    j_s(H_{\Lambda\cap\theta^\bot})\big]T^s \in K_0 (\mbox{\rm  Var}_\C)
    [[T]].
\label{auxiliarP}
 \end{equation}

 Notice that the Poincar\'e series $P(\Lambda)$ measures the classes
 in the Grothendieck ring of the jets of arcs in $H_\Lambda^*$ which
 are not jets of arcs in $H_{\Lambda \cap \theta^\bot}$, for any
 $0 \ne \theta \leq \s$, i.e., jets of arcs  with origin in $0$  which are not jets of arcs factoring through proper orbit
 closures of the toric variety $Z^\Lambda$.
It follows that:
\begin{Pro} \label{descompPgeom}
\[
P_{\mathrm {geom}}^{(Z^\Lambda,0)}(T)=\displaystyle\sum_{\theta\leq\s}P(\Lambda\cap\theta^\bot).
    \]
\end{Pro}
\begin{Exam} \label{P-sigma} The  series $P(\Lambda\cap\s^\bot)$ takes into account those
jets of arcs in $H_\Lambda$ which coincide with the jet of the constant arc. We have that
$P(\Lambda\cap\s^\bot)=\sum_{s\geq 0}[\{0\}]T^s=\sum_{s\geq 0}T^s$
hence
$
P(\Lambda\cap\s^\bot)=(1-T)^{-1}$.
\end{Exam}

\begin{Pro}
If $d =1$ and the multiplicity of the monomial curve $Z^\Lambda$
at the origin is equal to $m$ then the series $P(\Lambda)$ is
equal to $ {(\L-1)}{T^m}{(1-\L T)}^{-1}{(1-T^m)}^{-1}$.
\label{P_curvasmonomiales}
\end{Pro}

\begin{Exam} \label{P-codim1}  Let  $\Lambda$ be a semigroup, as in
Notation \ref{Lambda}, defining a toric variety of arbitrary dimension.
For any  $\theta \leq \s$ of codimension $1$ we denote by
$m_\theta$ the multiplicity of the monomial curve $(Z^{\Lambda \cap \theta^\bot},0)$.
Then we have
$P(\Lambda\cap\theta^\bot)= ({\L-1})  {T^{m_\theta}}  ({1-\L
T})^{-1} {(1-T^{m_\theta})^{-1}}$.
\end{Exam}

\begin{definition} Recall that $e_1 \dots, e_n$  denote the minimal
system of generators of the semigroup $\Lambda$.  The  {\em
$k^{th}$-logarithmic jacobian ideal} of $Z^\Lambda$ is the
monomial ideal $\J_k$ of $\C[\Lambda]$ corresponding to  the
following subset of $\Lambda$,
\begin{equation}       \label{j-k}
\{ e_{i_1} + \cdots + e_{i_k}  \;  \mid \; e_{i_1} \y \cdots \y
e_{i_k}  \ne 0,
  \mbox{ for }
1 \leq i_1 < \cdots <i_k \leq n  \}.
\end{equation}
\end{definition}
We abuse of notation by denoting also by $\J_k$ the set
(\ref{j-k}).
\begin{remark}
The motivation for this terminology is inspired by the Appendix in
\cite{LR} and  by the fact that these ideals are defined
geometrically in terms of differential forms on $Z^{\Lambda}$ with logarithmic
poles outside the torus (see Section \ref{tor}).
\end{remark}

\begin{Not} \label{jotak-bis}
 We denote
by $\Sigma_k $  (resp. by $\ord_{\mathcal{J}_k }$) the dual
subdivision of $\s$ (resp. the support function) of the Newton polyhedron of  the
$k^{th}$-logarithmic jacobian ideal  ${\mathcal J}_k $, for $k=1,
\dots, d$. The maps
\[
\left\{
\begin{array}{lcccccccl}
\f_1    &   :=  &   \mbox{\rm ord}_{\mathcal J_1  } &  \mbox{ and
} & \f_k  & :=  &   \mbox{\rm ord}_{\mathcal J_k  }  - \mbox{\rm
   ord}_{\mathcal J_{k-1}  } &  \mbox{ for }  k=2,\dots,d,
   \\
 \Psi_1 &:=& 0 & \mbox{ and } &  \Psi_{k} &:= & (k-1) \, \ord_{\J_k}-k \, \ord_{\J _{k-1}} & \mbox{ for }  k=2, \dots,
 d,
\end{array} \right.
\]
are piece-wise linear functions defined on the cone $\s$. If $\nu \in \s$ we put $\f_0 (\nu) := 0$
and $\f_{d+1} (\nu): = + \infty$.
If $\r \subset \s$ is a cone of dimension one, we denote by
$\nu_\r$ the generator of the semigroup $\r \cap N$. We define the
finite set:
\begin{equation} \label{BL}
    B(\Lambda) : = \left\{ (d,1) \right\}  \cup  \bigcup_{k=1}^d  \left \{ ({\Psi_k}(\nu_\r), { \f_{k}} (\nu_\r))
  \,   \mid \r \in\cup_{i=1}^k \Sigma_i^{(1)}, \mbox{ and } \stackrel{\circ}{\r} \cap
\stackrel{\circ}{\s} \ne \emptyset \mbox{ if } k<d
     \right\}    .
\end{equation}
 \end{Not}

\begin{remark}
Notice that the set $B(\Lambda)$ depends only on
 the Newton polyhedra (with integral structure) of the
logarithmic jacobian ideals  of $Z^{\Lambda}$.
In particular, we apply this observation to the
sets $B( \Lambda \cap \theta^\bot)$ for $\theta < \s$. For
$\theta = \s$ we convey that $B( \Lambda \cap \s^\bot) := \{ (0,
1) \}$.
\end{remark}

\begin{The}   \label{PLambdaRac}
The series $ P (\Lambda)$ is of the form
\[
P (\Lambda) =  Q_{\Lambda}  \prod_{(a, b) \in B(\Lambda) }
(1-\L^{a}T^{b})^{-1},  \mbox{ where }  Q_{\Lambda}  \in\Z[\L,T] \] is
determined by the lattice $M$ and the Newton polyhedra of the logarithmic jacobian ideals of $Z^\Lambda$.
\end{The}

\begin{Cor}   \label{P-geom}  With notations of Theorem  \ref{PLambdaRac}
the local geometric motivic Poincar\'e series of $(Z^\Lambda, 0)$,
\begin{equation}    \label{ratexp}
    P_{\mathrm {geom}}^{(Z^\Lambda, 0)}(T) = \sum_{\theta \leq \s}   Q_{\Lambda \cap \theta^\bot}
\prod_{(a, b) \in B(\Lambda \cap \theta^\bot) }
(1-\L^{a}T^{b})^{-1},
\end{equation}
is determined by the sequences of Newton polyhedra
of the logarithmic jacobian ideals of $Z^{\Lambda \cap \theta^\bot}$ and lattices
$M(\theta, \Lambda)$, for $\theta \leq \s$.
\end{Cor}

\begin{Cor}
Suppose that  the affine toric
variety $Z^\Lambda$ is normal.
If $\theta \leq \s$ we denote by
$\s_\theta^\vee$ the image of the cone $\s^\vee$ in
$(M_\theta)_\R$, where $M_\theta := M / \theta^\bot \cap M$ and by
$\Lambda (\theta)$ the semigroup $\Lambda (\theta) :=
(\s_\theta^\vee \cap M_\theta ) \times \Z_{\geq
0}^{\mathrm{codim}\theta}$. With this notation  we have
\[
 P_{\mathrm {geom}}^{Z^\Lambda}(T)=
 \sum_{\theta\leq\s}(\L-1)^{\mathrm{codim}\theta}
 P_{\mathrm {geom}}^{(Z^{\Lambda(\theta)},0)}(T).
\]
\label{P-locales}
\end{Cor}
\begin{remark}
See  also \cite{tesis} for a generalization of this result to the
class of affine toric varieties which are locally analytically
unibranched.
\end{remark}

We make more explicit the result for surfaces:

\begin{Cor}
Let $Z^\Lambda$ be an affine toric surface (case $d=2$ in Notation
\ref{Lambda}). We denote by  $\theta_1$ and $\theta_2$ the one
dimensional faces of the cone $\s$. The terms which appear in the denominator of
the rational expression (\ref{ratexp}) of $P_{\mathrm {geom}}^{(Z^\Lambda,0)} (T)$ are
$1-T,\ 1-\L T,\ 1-\L^2T$, $1-T^{m_{\theta_i}}$ and $1-\L^{\Psi_2
(\nu_\r)}T^{\phi_2 (\nu_\r)}$,
 where the integer $m_{\theta_i}$ is the multiplicity of the curve $Z^{\Lambda\cap\theta_i^\bot}$, for $i=1, 2$,
and  $\r $ runs through the rays of $ \Sigma_1 \cap  \Sigma_2$.
\label{poles-surfaces}
\end{Cor}

\begin{remark}
Suppose that  $\C$ denotes  the field of complex numbers. If $V$
is a variety the map $V \mapsto HD(V) \in \Z [u, v]$, where $HD(V)$
denotes the Hodge-Deligne polynomial,  factors through   $ K_0 (\mbox{\rm
  Var}_\C) $ inducing a ring morphism $HD:  K_0 (\mbox{\rm
  Var}_\C)   \rightarrow     \Z [u, v]$ which maps $\L \mapsto uv$ (see \cite{DL1}).
It follows that $\Z[\L] \cong \Z[X]$ where $X$ is an
indeterminate. By Corollary \ref{P-geom} the geometric motivic
Poincar\'e series of a toric singularity  is an element of
$\Z[\L](T)$, a ring in which the notion of the pole in  $T$ of a
non zero element is well defined since $\Z[\L]$ is an integral
domain.
\end{remark}
\begin{remark} \label{sur-rem}
 If $Z^{\Lambda}$ is a normal toric surface then  $T= 1$ and  $T = \L^{-1}$ are not poles of
$P_{\mathrm {geom}}^{(Z^\Lambda,0)}$  (see \cite{LR}). In general it
may happen
 that all candidate poles mentioned in the statement of Corollary \ref{poles-surfaces} are actual poles
 (see an example in Section \ref{toric-example}).
\end{remark}

We give a geometrical interpretation of the set of candidate poles
of the series $P_{\mathrm {geom}}^{(Z^\Lambda,0)} (T)$.

\begin{definition}
 For $1 \leq k \leq d$ we denote by $\p_{k}$ the composite of the
 normalization map $Z_{\s} \rightarrow Z^{\Lambda}$ with the toric
 modification of $Z_k \rightarrow Z_{\s}$ defined by the subdivision $\cap_{i=1}^k
 \Sigma_i$ of $\s$.
 \end{definition}

 The modification $\p_k$ is the minimal toric modification which factors through the normalization of $Z^{\Lambda}$
 and the normalized blowing up of $Z^\Lambda$  with center $\J_i$, for $i=1, \dots, k$.
 The rays $\r$ in the fan $\cap_{i=1}^k \Sigma_i$ correspond bijectively to
orbit closures of $Z_k$ which are of codimension one. If $\nu_\r$
is the generator of the semigroup $\r \cap N$ we denote by
$E_{\nu_\r}$ the irreducible component corresponding to $\r$. We denote by
$\mbox{val}_{{\nu_\r}}$ the divisorial valuation of the field of
fractions of $\C[\Lambda]$, which is associated to the divisor
$E_{\nu_\r}$. If $m \in M$ then we have that
\begin{equation} \label{val1}
\mbox{val}_{\nu_\r} (X^m) = \langle \nu_\r , m \rangle.
\end{equation}
We have that $\stackrel{\circ}{\r} \subset \stackrel{\circ}{\s}$
if and only if $E_{\nu_\r}$ is a codimension one irreducible
component of the exceptional fiber of $\p_k ^{-1} (0) $.
If $ 1 \leq i \leq k \leq d$, the pull-back $\p_k^* (\J_i)$ of
$\J_i$ by $\p_k$ defines a sheaf of locally principal monomial
ideals on the toric variety $Z_k$ and by (\ref{val1}) we deduce
\begin{equation} \label{val2}
 \mbox{val}_{\nu_\r} (\p_k^* (\J_i)) = \ord_{\J_i} (\nu_\r).
\end{equation}
\begin{Pro} \label{val3}
For $1 \leq k \leq d$,  \[ \mathcal{L}_k:= (\p_k^*
(\J_k))^{k-1} / (\p_k^* (\J_{k-1}))^{k} \mbox{ and } \mathcal{Q}_k
:= \p_k^* (\J_k)  / \p_k^* (\J_{k-1})\] are sheaves of locally
principal monomial ideals on $Z_k$ such that
\[
\begin{array}{lcl}
B(\Lambda)  & = &  \{ (d,1) \} \cup \bigcup_{k=1}^{d-1} \{
  ( \mbox{val}_{\nu_\r} ( \mathcal{L}_k ) , \mbox{val}_{\nu_\r}
  (\mathcal{Q}_k) ) \mid E_{\nu_\r} \subset \p_k ^{-1} (0) \}
  \\
 &  &  \cup   \{ ( \mbox{val}_{\nu_\r} ( \mathcal{L}_d ) , \mbox{val}_{\nu_\r}
  (\mathcal{Q}_d) ) \mid  \r \in \cap_{i=1}^d
  \Sigma_i, \dim \r = 1 \}.
\end{array}
\]
\end{Pro}

\begin{remark}
If $(S,0)$ is equidimensional of dimension $d$ then the term $1-
\L^d T$ appears always in the denominator of the rational form of  $ P_{\mathrm
{geom}}^{(S,0)} (T)$. This is consequence of Theorem 7.1 of
\cite{DL1}.
\end{remark}

\section{Combinatorial convexity properties of Newton
polyhedra of $\J_k$} \label{conv-newton}

We study the combinatorial convexity properties of the support
functions of the Newton polyhedra of the monomial ideals $\J_k
\subset \C[ \Lambda]$, for $\Lambda$ as in Notation \ref{Lambda}.

 If $\nu \in \s$ then  the relation   $\leq_{\nu}$ defined by
\begin{equation}  \label{partial-nu}
v \leq_{\nu} v' \Leftrightarrow  \langle \nu, v \rangle \leq
\langle \nu, v' \rangle,
\end{equation}
is a preorder on the set $\s^\vee \cap M$. We give an
algorithm to determine a vector $w_k \in \J_k  $ such that
$\ord_{\J_k}    (\nu) = \langle \nu, w_k \rangle$.
\begin{Lem}    \label{alg-jotak}
Let $\nu$ be an element of $\stackrel{\circ}{\s} \cap {N} $ such that
\begin{equation}
e_1 \leq_\nu e_2 \leq_\nu \cdots  \leq_\nu e_n,
\label{ordenvectores}
\end{equation}
for the preorder $\leq_\nu$ defined by (\ref{partial-nu}).
Define the sequence $i_1 < i_2 < \cdots < i_d \leq n $ of $1,
\dots, n$ in the following inductive form: set $i_1:=1$, suppose
that $i_2, \dots, i_k$ have already been defined and set:
\begin{equation} \label{pot}
i_{k+1} : = \min \{ 1 \leq i \leq n  \ |\ e_{i_1} \y \cdots \y
e_{i_{k}} \y e_i \ne 0 \}.
\end{equation}
Set $w_k :=e_{i_1} + \cdots + e_{i_k}$ for $k=1, \dots,d $. Then
we have:
\begin{equation} \label{prepa}
\ord_{\J_k} (\nu) = \langle \nu, w_k \rangle \quad \mbox{ and }
\quad \f_k    (\nu) = \langle \nu, e_{i_k} \rangle.
\end{equation}
\end{Lem}
{\em Proof.} We deduce from (\ref{pot}) that:
\begin{equation} \label{put}
\langle \nu, e_{i_{k+1}} \rangle = \min  \{ \langle \nu, e_{i}
\rangle \ |\ 1\leq i\leq n, \exists 1\leq {j_1}<\cdots <j_k\leq n,\ e_{j_1}
\y \cdots \y e_{j_{k}} \y e_i \ne 0 \}.
\end{equation}

The statement is obvious for $k=1$.  Suppose the result for $1 <k
< d$. We have then that:
\[
\ord_{\J_{k+1}}  (\nu ) \stackrel{\mbox{\rm \small Def.}}{\leq}
\langle \nu, e_{i_1} + \cdots + e_{i_{k+1}} \rangle
\stackrel{\mbox{\rm \small by
    induction}}{=}
\ord_{\J_k} (\nu) + \langle \nu, e_{i_{k+1}} \rangle.
\]
The other inequality follows from Formula (\ref{put}) since
\[
\begin{array}{lcl}
\ord_{\J_{k+1}}  (\nu ) & =  &  \min \{ \langle \nu, e_{j_1} +
\cdots + e_{j_{k+1}} \rangle \}_{0 \ne  e_{j_1} \y \cdots \y
e_{j_{k+1}} }
\\
& \geq & \min \{  \langle \nu, e_{j_1} + \cdots + e_{j_{k}}
\rangle \}_{0 \ne   e_{j_1} \y \cdots \y e_{j_{k}} }  + \min
\langle \nu, e_{j_{k+1}} \rangle
\\
& = & \ord \J_{k}  (\nu ) + \langle \nu, e_{i_{k+1}} \rangle.
\end{array}
\]\hfill $\ {\Box}$

\begin{Pro}  \label{algo}
For every $\nu$ in $\stackrel\circ\s\cap N$ there exist $1\leq i_1, \dots, i_d \leq n$
such that $\f_k (\nu) = \langle \nu, e_{i_k} \rangle $, $\sum_{r=1}^k e_{i_r}   \in \J_k$ and
 $\ord_{\J_k}    (\nu) =  \langle \nu, \sum_{r=1}^k e_{i_r}    \rangle$, for $k=1, \dots, d$.
\end{Pro}

{\em Proof.} It follows immediately from Lemma \ref{alg-jotak}. \hfill $\Box$

\begin{Cor} \label{c1-bis}
If $\nu \in \stackrel{\circ}{\s} \cap {N} $ we have that:
\begin{equation} \label{menor-bis}
0 = \f_0 (\nu)  < \f_1 (\nu) \leq \f_2 (\nu) \leq \cdots \leq \f_d (\nu) < \f_{d+1} (\nu) = + \infty.
\end{equation}
\end{Cor}

\begin{definition} \label{ak}
For $0 \leq k \leq d$ we set $
A_k := \{ (\nu,s)  \, \mid  \, \nu \in
\stackrel{\circ}{\s} \cap N,  \, \,  \f_k (\nu)  \leq    s   < \f_{k+1} (\nu)  \}$.
\end{definition}

\begin{Pro}
 \label{c2-bis}  The sets  $A_0, \dots, A_d$  define a partition of
$(\stackrel{\circ}{\s} \cap N) \times \Z_{> 0}$.
\end{Pro}

{\em Proof.} It follows from Corollary \ref{c1-bis}. \hfill $\Box$

\begin{definition} \label{ele}
If $(\nu,s) \in  (\stackrel{\circ}{\s} \cap N) \times \Z_{> 0}$ we denote by $\ell_\nu^{s}$ the linear subspace of
$M_\Q$ given by
\[\ell_\nu^{s}  :=  \mbox{\rm span}_\Q \{ e_i \mid 1 \leq i \leq n \mbox{ and }
\langle \nu, e_i \rangle \leq s \}.\]
\end{definition}

\begin{Lem} \label{crucial-bis}
Let  $(\nu,s) $ belong to  $A_k$ for some $1 \leq k \leq d$.  Let $w_k
\in {\mathcal J}_{k}  $ verify that  $ \mbox{\rm ord}_{{\mathcal J}_{k}  } (\nu)  = \langle
\nu, w_k \rangle $. If $w_k = e_{j_1} + \cdots + e_{j_k}$
is an expansion as a sum of  $k$ linearly independent vectors in $\{e_1, \dots, e_n \}$
then  $\{ e_{j_1}, \dots , e_{j_k}
\}$  is a basis of  the vector space
$\ell_\nu^{s}$.
If $
e_i \in  \ell_\nu^{s} $  and
$e_i = \sum_{r=1}^k \a_r e_{j_r}$
then $\a_r \ne 0$  implies that $\langle \nu, e_{j_r} \rangle \leq
\langle \nu, e_i \rangle$, for $r= 1,\dots,k$.
\end{Lem}
{\em Proof.}
Let  $e_{i_1}, \dots, e_{i_{d}}$  be the vectors defined  by Proposition \ref{algo}.
We set $w_r' := e_{i_1} + \cdots + e_{i_r}$
for $r= 1, \dots, d$.   By Proposition \ref{algo}  and   Corollary \ref{c1-bis}
if $1 \leq r \leq k$ then we deduce
$\f_r (\nu ) = \langle \nu,e_{i_r} \rangle \leq s$.
               This implies that $\mbox{\rm span}_\Q \{
e_{i_1},  \dots,  e_{i_{k}}  \}  \subset \ell_\nu^s$.
If for some $1 \leq i \leq n$,  $\langle
\nu, e_i \rangle \leq s$ and $e_i \notin    \mbox{\rm span}_\Q \{
e_{i_1},  \dots,  e_{i_{k}}  \} $ then the
vector $\bar{w}_{k+1}: = e_{i_1} + \dots + e_{i_{k}} + e_i $
belongs to $\J_{k+1} $  hence
$
\ord_{\J_{k+1}} ( \nu) \leq \langle \nu, \bar{w}_{k+1} \rangle  =
\ord_{\J_{k}} ( \nu)
 + \langle \nu, e_i \rangle$.
This implies that $\f_{k+1} (\nu) \leq \langle \nu, e_i \rangle
\leq s$, a contradiction with the fact that $(\nu,s) \in {A}_k$. We have shown
that $\ell_\nu^s =  \mbox{\rm span}_\Q \{
e_{i_1},  \dots,  e_{i_{k}}  \} $.

Suppose  that there exists a vector
 $w_k = e_{j_1} + \cdots + e_{j_k}   \in   \J_k$ such that:
\begin{equation} \label{inek}
\ord_{\J_k} (\nu) = \langle \nu, w_k ' \rangle  = \langle \nu, w_k
\rangle \quad \mbox{\rm and }  \quad \mbox{\rm span}_\Q \{
e_{i_1},  \dots,  e_{i_{k}}  \} \ne \mbox{\rm span}_\Q \{ e_{j_1},
\dots,  e_{j_{k}}  \}.
\end{equation}
Then there exists $1 \leq k_0 \leq k$ such that $e_{j_{k_0}}
\notin \mbox{\rm span}_\Q \{   e_{i_1},  \dots, e_{i_{k}}  \}$. If
$\langle \nu, e_{j_{k_0}} \rangle  < \langle \nu, e_{i_k} \rangle$
then the vector $\hat{w}_k :=  e_{i_1} + \cdots + e_{i_{k-1}} +
e_{j_{k_0}}$ belongs to $\J_k$ and $\langle \nu, \hat{w}_k \rangle
< \ord_{\J_k} (\nu)$ by (\ref{inek}), a contradiction.  If
$\langle\nu,e_{i_k}\rangle<\langle\nu,e_{j_{k_0}}\rangle$, then
the vector $w_k'-e_{j_{k_0}}$ belongs to $\J_{k-1}$ and we have
$\langle\nu,w_k'-e_{j_{k_0}}\rangle<\langle\nu,w_{k-1}' \rangle$,
which contradicts the formula
$\ord_{\J_{k-1}}(\nu)=\langle\nu,w_{k-1}'\rangle$. Thus, the equality
\begin{equation} \label{rbq}
\langle \nu, e_{j_{k_0}} \rangle  =  \langle \nu, e_{i_k} \rangle
\end{equation}  holds. The equality $\f_k (\nu) = \f_{k+1} (\nu)$ follows
from (\ref{rbq}) and  the inequalities:
\[
\langle \nu, e_{i_k} \rangle \stackrel{\mbox{\rm \small
    (\ref{prepa})}}{=} \f_k (\nu)  \stackrel{\mbox{\rm \small
    (\ref{menor-bis})}}{\leq}  \f_{k+1} (\nu) \stackrel{\mbox{\rm \small
   Def. }}{=}  \ord_{\J_{k+1}}
    (\nu) - \ord_{\J_{k}} (\nu),
\]
\[
\ord_{\J_{k+1}} (\nu) \stackrel{\mbox{\rm \small
   Def. }}{\leq}  \langle \nu, e_{i_1} + \cdots + e_{i_k}  +
e_{j_{k_0}} \rangle  \stackrel{\mbox{\rm \small
    (\ref{prepa})}}{=}   \ord_{\J_{k}} (\nu)  + \langle \nu, e_{j_{k_0}} \rangle.
\]

Finally, if  we have an
expansion $ e_i = \sum_{r=1}^k \a_r e_{j_r}$ with  $\a_{r_0} \ne 0$
and $\langle \nu, e_i \rangle < \langle \nu, e_{j_{r_0}} \rangle$
then the vector $ \tilde{w}_k: = e_i + \sum_{r=1, r\ne r_0}^k
e_{j_r} $ belongs to $\J_k$. The inequality
$
\ord_{\J_k} (\nu) = \langle \nu, e_{j_1} + \cdots + e_{j_k}
\rangle > \langle \nu , \tilde{w}_k\rangle,
$
is a contradiction
with the definition of the support function. \hfill $\ {\Box}$

\section{The universal family of arcs in the torus}
\label{universal}

We describe the universal family of arcs in the torus $T_N$ of a
rank $d$ lattice $N$ and some properties of its functions which
are useful to deal with jets of arcs in toric varieties.
In this section we use that the characteristic of the base
field $\C$ is zero.

Let us fix a basis $m_1, \dots, m_d $ of $M$. We set
\[ \mathcal{A} :=  \C[ c(m_i)^{\pm 1} ] \otimes_\C \C[ u_j(m_i)]_{i=1,\dots,
d}^{j\geq 1},  \]
 where   $\{ c(m_1),\dots, c(m_d) \} \cup \{ u_j(m_i) \}_{i=1,\dots,
d}^{j\geq 1}$ are algebraically independent over $\C$. Then there
is one homomorphism of semigroups
 $h^*: (M, +) \rightarrow   (\mathcal{A} [[t]]^*, \cdot)$ such that
$m_i \mapsto c(m_i) ( 1 + \sum_{j \geq 1} u_j (m_i) t^j)$ for
$i=1, \dots, d$. We have that the image of $m \in M$ by this
homomorphism is to $c(m) u(m)$ where $u(m) \in \mathcal{A}
[[t]]$ is a series of the form $ u(m) = 1 + \sum_{j\geq 1} u_j( m)
t^j$. Notice that if ${m}, {m}' \in M$ then we have $ c({m} +
{m}') = c({m}) c({m}')$ and $ u(m + m') = u(m) u(m')$. By the
description of  Section \ref{intro-arcs} we check that    $\mathcal{A} = \mathcal{A}_{T_N}$
and the map $h: \mbox{\rm Spec}  \mathcal{A} [[t]]  \rightarrow
T_N$ corresponding to $h^*$,  is the universal family of arcs in
the torus.
The following two lemmas show some relations among the elements
$u_i ({m}) \in \mathcal{A}_{T_N}$, when we vary $i$ and ${m} \in
M$, in terms of linear dependency relations among the ${m} \in M$.

\begin{Lem} \label{basis}
Let $ {m}_1, \dots, {m}_k$ be a set of linearly independent
vectors of $M$. If $0 \ne {m} = \sum_{j=1}^k a_j {m}_j $  with $a_j \in
\Z$ then we have:
\begin{equation} \label{lam}
u_i ({m}) =  \sum_{r=1}^k a_r \, u_i ({m}_r) + R^{({m})}_{i-1}
(u_j({m}_r ))^{j=1, \dots, i-1}_{r=1, \dots, k},
\end{equation}
for all $i \geq 1$,  where $R^{({m})}_{i-1}$ is a
quasi-homogeneous polynomial of weight $i$, with rational
coefficients in $(u_j({m}_r ))^{j=1, \dots, i-1}_{r=1, \dots, k}$,
where  $u_j({m}_r)$ is given weight equal to $j$, for $r=1, \dots,
k$ and $j=1, \dots, i-1$.
\end{Lem}
{\em Proof.} We have that $
 {u} ( {m} ) =  \prod_{j=1}^k {u} ({m}_j)^{a_j} =    \prod_{j=1}^k
 \left(  \sum_{i \geq 0} u_i ({m}_j) t^i \right)^{a_j}
$. Remark that if $\phi= \sum u_i t^i \in \C[[t]]$ is a series
with constant term equal to one and if $n \in \Z$ then the series
$\phi^n$ is of the form: $ \phi^n = \sum P^{(n)}_i t^i $ where
$P^{(n)}_i = n  u_i + R^{(n)}_i (u_1, \dots, u_{i-1} )$ is a
quasi-homogeneous polynomial in $u_1, \dots, u_{i}$, where  $u_l$
has weight equal to $l$; notice that the  coefficient $n$ of $u_i$ 
does not vanish since $\C$ is a field of characteristic zero. 
We use this observation to compute the
expansion of $ \prod_j  \left(  \sum_{i \geq 0} u_i ({m}_j)
t^i\right)^{a_j}$ as a series in $t$. \hfill $\ {\Box}$

\begin{Lem} \label{vector}
If ${m'}_1, \dots, {m'}_k$ and ${m}_1, \dots, {m}_k$ are linearly
independent vectors in the lattice $M$ spanning the same linear
subspace $\ell$ of $M_{\Q}$ then for any $s \geq 1$  we have the
equality of $\Q$-algebras:
\begin{equation} \label{kalg}
\Q [ u_1( {m'}_j ), \dots, u_s ({m'}_j ) ]_{j=1}^k =\Q [ u_1(
{m}_j ), \dots, u_s ({m}_j ) ]_{j=1}^k.
\end{equation}
In particular, if ${m} = \sum_{j=1}^k a_j m_j$, with $a_j \in \Q$
then  $u_i ({m}) $ belongs to the $\Q$-algebra (\ref{kalg})
for $i=1,\dots, s$.
\end{Lem}
{\em Proof.} It is sufficient to prove it in the case that ${m}_1,
\dots, {m}_k$ are a basis of the rank $k$ lattice $\ell \cap M$.
We show the result by induction on $s$.  Since ${m}_1, \dots,
{m}_k$ form a basis of $\ell \cap M$ and ${m'}_r \in \ell \cap M$
we have expansions:
$
{m'}_r= a_{r,1} {m}_1 + \cdots + a_{r, k} {m}_k$
with $a_{r,j} \in \Z $, for
$ j = 1, \dots, k$,
and
$r = 1, \dots, k$.
Since ${m'}_1, \dots, {m'}_k$ are linearly independent we have
expansions:
$
{m}_r= b_{r,1} {m'}_1 + \cdots + b_{r, k} {m'}_k $
$b_{r,j} \in \Q $,
for $ j = 1, \dots, k$ and
$r = 1, \dots, k$.
For $s=1$  the term $R_0^{(m)} $ appearing in formula
(\ref{lam}) is equal to zero thus by Lemma \ref{basis} we obtain
that
  $u_1 ({m'}_r) = a_{r,1} u_1 ({m}_1)  + \cdots + a_{r, k} u_1({m}_k)  $
and
   $u_1( {m}_r) = b_{r,1} u_1( {m'}_1)  + \cdots + b_{r, k} u_1(
{m'}_k )$, for $ r =1, \dots, k$.
Using the induction hypothesis for all $1 \leq s' <s$ and the
triangular  form of formula
 (\ref{lam}) for $u_s ({m'}_r)$
we deduce that $u_s ({m}_r)$ is of the form:
\[
u_s( {m}_r) = b_{r,1} u_s( {m'}_1)  + \cdots + b_{r, k} u_s(
{m'}_k ) + P_{r, s},
\]
where $P_{r, s}$ belongs to $\Q [ u_1 ( {m'}_j), \dots  u_{s-1} (
{m'}_j) ]_{j=1}^k$. \hfill $\ {\Box}$

\begin{Pro} \label{alg-ind}
If the vectors $ {m}_1, \dots, {m}_k$ in $M$ are linearly
independent  then the following elements of $\mathcal{A}_{ T_N}$
are algebraically independent over $\C$:
\begin{equation} \label{lla}
c({m}_1), \dots, c({m}_k), \mbox{ and }
 u_i({m}_1), \dots,
u_i({m}_k), \quad \forall i \geq 1
\end{equation}
\end{Pro}

\section{The image in the Grothendieck ring of the jets of the
orbits} \label{tor-jet}

Let $\nu$ belong to  the set $\stackrel{\circ}{\s} \cap N$. We consider the orbit $H^*_\nu$ of the
action of the arc space of the torus on $H^*_\Lambda$.
The universal family of arcs in the
torus parametrizes the arcs in $H^*_\nu$ by  the morphism
$\Psi_\nu: \mbox{{\rm Spec}} \mathcal{A}_{T_N} [[ t]] \rightarrow
Z^\Lambda$ given by:
\[
X^{e_i} \circ \Psi_\nu = t^{\langle \nu, e_i \rangle }  c( e_i) u
( e_i), \mbox{ for } i = 1, \dots, n.
\]
Recall that $e_1, \dots, e_n$ is the minimal system of generators
of $\Lambda$ (cf. Notations \ref{Lambda}).

We  prove that the set $j_s( H^*_{\nu})$ is a locally closed
subset of $j_s( \A^n_\C)_0 \cong \A^{ns}_\C$ and we determine its class in the
Grothendieck ring of varieties.

\begin{The}\label{key-bis}
If $(\nu, s) \in A_k$ for some $0 \leq k \leq d$, then
    the jet space $j_s (
H^*_{\nu})$          is  a
locally closed subset of $H_s( Z^\Lambda )_0$ isomorphic to $\{0 \}$ if $k =0$ or to
  $(\C^*)^{k}
\times \A_\C^{k s  - \mbox{\rm \small ord}_{\mathcal J  _{k}}
(\nu)}$ if $k >0$.
\end{The}
{\em Proof.}
 If
$h \in H^*_\nu $   the equality  $ \ord_t
 ( X^{e_i} \circ h ) = \langle \nu, e_i \rangle $ holds for $1 \leq i \leq n$.
By Definition      \ref{ele} those vectors $e_i$ such that $j_s( X^{e_i} \circ h) \ne 0$
span the     $\Q$-vector space    $\ell_\nu^s$ since $ \langle \nu, e_i \rangle \leq s$.
 If $k=0$ this vector space is empty, the jet space $j_s (
H^*_{\nu})$ consists of the constant $0$-jet and the conclusion follows.

Suppose then that $k >0$. 
We denote by 
$\mathcal{O}_\nu^s$     (resp. by $C_\nu^s$)
the $\C$-algebra of $\mathcal{A}_{T_N}$ generated by:
 \begin{equation} \label{gener}
 c(e_i)^{\pm 1}, u_1 (e_i), \dots , u_{s - \langle \nu, e_i \rangle} (e_i)
 \mbox{ for those } 1\leq i\leq n  \mbox{ such that }   \langle \nu,
 e_i \rangle \leq s, 
 \end{equation}
\[
\mbox{(resp.   }    \, c(e_i)^{\pm 1}                 \mbox{ for those } 1\leq i\leq n  \mbox{ such that }   \langle \nu,
 e_i \rangle \leq s \mbox{ \ ).} 
\]
By Proposition \ref{algo} the vector $\nu$ determines 
integers $1 \leq i_1, \dots, i_k \leq n$ such that  
 $\ord_{\mathcal J  _{k}} (\nu) = \sum_{r=1}^k \langle \nu, e_{i_r} \rangle$.

\par
\noindent
\textbf{Assertion.} \textit{We have the following properties}:
\begin{enumerate}
 \item[(i)] \textit{ Denote by $\underline{U}$ the variables $(U_1, \dots, 
U_{k s - \mbox{\rm \small ord}_{\mathcal J_{k} (\nu)}})$.  For any $1 \leq i \leq n$ and $l$ such that 
$1 \leq l \leq s - \langle \nu, e_i \rangle$ 
there exists a polynomial $P_{l, i} \in  \Q [ \underline{U} ]$ such that    }
\[
u_l (e_i) = P_{l, i} \left(u_1 (e_{i_1}), \dots, u_{s - \langle \nu, e_{i_1} \rangle} (e_{i_1}), \dots, 
u_1 (e_{i_k}), \dots, u_{s - \langle \nu, e_{i_k} \rangle} (e_{i_k}) \right). 
\]

\item[(ii)]      \textit{ The ring $\mathcal{O}_\nu^s$  is generated as a $C_\nu^s$-algebra by
\begin{equation} \label{pre-gen}
u_1 (e_{i_r}), \dots, u_{s - \langle \nu, e_{i_r} \rangle }
(e_{i_r}), \mbox{ for } r=1, \dots, k.
\end{equation}
  }
\item[(iii)] \textit{  The lattice $M_\nu^s$ spanned by 
$\{ e_i \mid 1 \leq i \leq n, \ \langle \nu, e_i \rangle \leq s \} $ is of rank $k$. 
The map $C_\nu^s \to \C [M_\nu^s]$ given by $c(e_i) \mapsto X^{e_i}$ is an isomorphism.  }

\item[(iv)] \textit{ The variety $\mathrm{Spec} \mathcal{O}_\nu^s$ is isomorphic to    $(\C^*)^{k}
\times \A_\C^{k s  - \mbox{\rm \small ord}_{\mathcal J  _{k}} 
(\nu)}$             } .
\end{enumerate}
{\em Proof of the Assertion.}
By Lemma  \ref{crucial-bis} the vectors  $e_{i_1},\ldots,e_{i_k}$  define a basis of $\ell_\nu^s$.
If $e_i$ is a vector in  $\ell_\nu^s$ with $\langle \nu, e_i \rangle \leq s$, and 
  $e_i = \sum_{r=1}^k \a_r e_{i_r} \in \ell_\nu^s$, then  $\a_{r} \ne 0$ implies that
$\langle \nu, e_i \rangle \geq \langle
\nu, e_{i_r} \rangle$.
We deduce from Lemma \ref{vector} that the elements $u_1
(e_i), \dots, u_{s - \langle \nu, e_i \rangle} (e_i) $ belong to
the $\Q$-algebra generated by (\ref{pre-gen}).   This implies that (i) and (ii) hold.

By Lemma \ref{alg-ind} and the definitions the ring $C_\nu^s$  is isomorphic to 
the $\C$-algebra of the lattice $M_\nu^s$. 
The assertion (iii) follows,  since by definition  $M_\nu^s$ is a sublattice of finite index of the rank $k$ 
lattice $\ell^s_\nu \cap M$. 

Finally, (iv) follows from these observations and Lemma \ref{alg-ind}. This ends the proof of the Assertion.   
\par

By the Assertion the morphism 
\[
\psi: \mathrm{Spec} \ C_\nu^s \otimes_\C  \C [\underline{U}] \longrightarrow j_s (\A^n_\C)_0, 
\]
given by
\[
 x_i (t) := \left\{ 
\begin{array}{ccl}
 c(e_i) \, t^{ \langle \nu, e_i \rangle } \,  \left( 1 + P_{1, i} (\underline{U}) t + \cdots P_{s -\langle \nu, e_i \rangle, i} (\underline{U})
t^{s - \langle \nu, e_i \rangle}  \right)  & \mbox{ if } &     \langle \nu, e_i \rangle \leq s, 
 \\
0 &        \mbox{ if } &   \langle \nu, e_i \rangle > s, 
\end{array}  \right. 
\]
for $1 \leq i \leq n$, is an inmersion. The image  $\mathrm{Im} (\psi)$ of $\psi$ is a locally closed subset. 

Finally,  if $h$ belongs to $H^*_\nu$, then $j_s (h)$ belongs to  $\mathrm{Im} (\psi)$ by the Assertion. 
Conversely, if $\xi \in \mathrm{Im} (\psi)$  we define an arc $h  \in H^*_\nu$ such that $j_s (h) = \xi$
 by specialization
from the universal family. 
First, for $1 \leq i \leq n$ and $l \geq 1$ we set 
$u_{l} (e_i) = 0$ if $l  > s - \langle \nu, e_i \rangle$.
By the Assertion the coefficients    $u_{l} (e_i)$ associated to $h$,  for     
 $ 1 \leq s - \langle \nu, e_i \rangle \leq  l$, are complex numbers 
determined by $\xi$. In order to complete the definition of  $h$ we have to 
give values for the coefficients $c(e_i)$ corresponding to $h$. We have 
an injection of $\C$-algebras $\C[M_\nu^s] \subset \C [M] =    \C [c(e_i)^{\pm 1}]_{i=1}^n$ which corresponds to 
a surjective map of torus. The inicial coefficients associated to $\xi$ define a closed point $p(\xi)$ of the torus 
$\mathrm{Spec} \C[M_\nu^s]$. Any closed point in the fiber of $p(\xi)$ by this map provides suitable initial 
coefficients $c(e_i) \in \C^*$, in such a way that the resulting arc $h$ verifies that $j_s (h) = \xi$.

Notice that $\mathcal{O}_\nu^s$ is the coordinate ring of the locally closed subset $j_s (H^*_\nu)$. 
  \hfill $\ {\Box}$

\section{Description of the series $P(\Lambda)$}

\label{GenptOrbits}

We describe the coefficients of the auxiliary
series $P(\Lambda)$.
We study in which cases  the intersections $j_s (H^*_\nu) \cap j_s
(H^*_{\nu'})$ and  $j_s(H_{\Lambda, \nu}^*)\cap
j_s(H_{\Lambda\cap\theta^\bot})$ are  non-empty, for $(\nu,s), (\nu',s) \in A_k$ and $\theta \leq \s$.

\begin{definition} \label{eq}
Define an equivalence relation in the set ${A}  _k $  for any $1
\leq k \leq d$:
\begin{equation} \label{equivalence-bis}
(\nu,s)  \sim (\nu', s')  \in {A}_k   \Leftrightarrow \left\{
\begin{array}{c}
 s = s ' ,\  \nu  \mbox{ and  }  \nu' \mbox{ define the same face of }
{\mathcal N}( \mathcal{J}_j)
\\
\mbox{ and  } {\mbox{\rm ord}}_{{\mathcal J}_j} (\nu) = {\mbox{\rm
ord}}_{{\mathcal J}_j} (\nu')  \mbox{, for } 1 \leq j \leq k.
\end{array}
 \right.
\end{equation}
We denote by $[(\nu,s)]$ the equivalence class of $(\nu, s)$ in
$A_k$ by this relation.
\end{definition}

\begin{remarks} \label{remk} $\,$
\begin{enumerate}
\item[(i)] For any fixed integer $s_0 >0$ the set $\{ [ (\nu, s_0)
] \mid
 (\nu, s_0) \in A_k \} $ is finite for $1\leq k\leq d$.

\item[(ii)] If $k=d$ the equivalence relation defined in the set
$A_d$ is the equality.
\end{enumerate}
\end{remarks}

\begin{Pro} \label{Prop-equiv}
If $(\nu,s), (\nu',s) \in A  _k$ the following relations are
equivalent:
\begin{enumerate}
\item[(i)] $(\nu,s) \sim (\nu',s)$,

\item[(ii)]  $\ell_\nu^{s} =\ell_{\nu'}^{s} $ and $\nu_{|
\ell^{s}_\nu }=\nu'_{| \ell^{s}_{\nu'} },$

\item[(iii)] $j_s(H_{\nu}^*)=j_s(H_{ \nu'}^*)$,

\item[(iv)] $j_s(H_\nu^*)\cap j_s(H_{\nu'}^*)\neq\emptyset.$
\end{enumerate}

\end{Pro}
{\em Proof.}
The condition $(\nu,s)\sim(\nu',s)$  implies that
$\ell_\nu^s=\ell_{\nu'}^s$ by Lemma \ref{crucial-bis}.
 The condition $\mbox{\rm ord}_{\J_j}(\nu)=\mbox{\rm
ord}_{\J_j}(\nu')$ for $j=1,\ldots,k$ is equivalent to
$\nu_{|_{\ell_\nu^s}}=\nu'_{|_{\ell_{\nu'}^s}}$ by Lemma
\ref{alg-jotak}. It follows that the conditions (i) and (ii) are
equivalent.

If (ii) holds then the basis $e_{j_1},\dots, e_{j_k}$ of the
vector space introduced  in Lemma \ref{alg-jotak}, coincides
for the vectors $\nu$ and $\nu'$ and  $\langle \nu, e_{j_r}
\rangle = \langle \nu', e_{j_r} \rangle$, for $r=1, \dots, k$.
This implies that the inmersion $\psi$ of $j_s
(H^*_\nu)$ defined for $(\nu, s)$ in the proof of Theorem
\ref{key-bis} is the same map as the one defined for $(\nu',s)$,
hence $j_s (H^*_\nu) = j_s (H^*_{\nu'})$.

Suppose that (iii) or (iv) holds. If  $h \in H^*_\nu$, $h' \in
H^*_{\nu'}$ verify that   $0 \ne j_s( X^{e_i}  \circ h)  = j_s(
X^{e_i} \circ h') $ for some  $1 \leq i \leq n$, then $X^{e_i}
\circ h$ and $X^{e_i} \circ h'$ have the same order $\langle \nu,
e_i \rangle = \langle \nu', e_i\rangle \leq s$ and those vectors $e_i$
generate the linear subspace $\ell_\nu^s = \ell_{\nu'}^s$,  hence (ii) holds. \hfill $\ {\Box}$

\begin{Not} \label{auxi}
The cone $\hat{\s}: = \s \times \R_{\geq 0}$ is rational for the
lattice $\hat{N}: = N \times \Z$.
\begin{enumerate}
\item[(i)]  If $\t \subset \s$ and $1 \leq k \leq d $ we set
$\t(k) := \{ (\nu, s) \mid \nu \in \stackrel{\circ}{\t} \cap
\stackrel{\circ}{\s} \mbox{ and } \f_k( \nu) \leq s < \f_{k+1}
(\nu) \}$.

\item[(ii)] If $\t \in \cap _{i=1}^k \Sigma_i$ then we set $A_{k,
\t}: = \t(k) \cap \hat{N}$.
\end{enumerate}
\end{Not}

\begin{remark} \label{pwl}
If $\t$ is a cone contained in a cone of the fan $\cap_{i=1}^k
\Sigma_i$ and if  $\t(k) \ne \emptyset$, then the closure of
$\t(k)$ in $\hat{\s}$ is  a convex polyhedral cone, rational for the lattice $\hat{N}$
(since in this case
the functions $\ord_{\J_1}, \dots. \ord_{\J_k}$ , hence also $\f_1, \dots, \f_k$ , are linear on $\t$  and the function $\ord_{\J_{k+1}}$, hence also $\f_{k+1}$, is piece-wise linear and  convex on
$\t$). In particular, the set
$A_{k, \t}$ may be empty, for instance, if $\t$
is contained in the boundary of $\s$ or if for all $\nu$ in the
interior of $\t$ we have that $\f_k(\nu) = \f_{k+1} (\nu)$.
\end{remark}

\begin{remark} \label{86}
 We deduce the following:
\begin{enumerate}
\item[(i)]
 $A_{k} =
{\sqcup_{\t \in \cap _{i=1}^k \Sigma_i}} A_{k,
\t}$ for $1 \leq k \leq d$ .

\item[(ii)] The vectors  $(\nu,s), (\nu',s) \in A_k$ are
equivalent by the relation $\sim$ in (\ref{equivalence-bis}) if and
only if there exists a cone $\t \in  \cap _{i=1}^k \Sigma_i$ such
that $\nu$ and $\nu'$ belong to the relative interior of $\t$ and
 $\f_i (\nu) = \f_i (\nu')$ for $i= 1, \dots, k$.

\item[(iii)]  It follows that $A_k /_{\sim} =
{\sqcup_{\t \in \cap _{i=1}^k \Sigma_i}} A_{k, \t}
/_{\sim}$, where  $A_{k, \t} /_{\sim}$ is the set of equivalent
classes of elements in the set $A_{k, \t}$ by the relation (\ref{equivalence-bis}).
\end{enumerate}
\end{remark}

\begin{Pro}
If $\nu \in \stackrel{\circ}{\s}$, $s \geq 1$  and $\theta \leq
\s$ then the following relations are equivalent:
\begin{enumerate}

 \item[(i)] $j_s(H_{\Lambda,
\nu}^*)\cap j_s(H_{\Lambda\cap\theta^\bot}) \ne \emptyset$,

\item[(ii)] $j_s(H_{\Lambda, \nu}^*)\subset
j_s(H_{\Lambda\cap\theta^\bot})$,

\item[(iii)] $\ell_\nu^s \subset\theta^{\bot}$.
\end{enumerate}
\label{thetaproposition}
\end{Pro}
\noindent{\em  Proof.} Suppose that (i) holds. Then there is an arc $h
\in H^*_\nu$ whose $s$-jet belongs to
$j_s(H_{\Lambda\cap\theta^\bot})$. If the truncation
 $j_s (X^{e_i} \circ h)$ does not vanish then the
vector $e_i$ belongs to $ \Lambda \cap \theta^\bot$.
    By Definition \ref{ele},
those vectors $e_i$ for which       $j_s (X^{e_i} \circ h) \ne 0$
span the linear subspace  $\ell_\nu^s$. This proves the inclusion
(iii).

Assume that (iii) holds.  Let $h \in  H_{\Lambda,\nu}^*$. Define
an arc $h' \in H_{\Lambda \cap \theta^\bot}$ by the semigroup
homomorphism $ \Lambda \cap \theta^\bot \rightarrow \C[[t]]$,
given by $e \mapsto X^e \circ h$, for $e \in   \Lambda \cap
\theta^\bot $. We have that $h' \in H^*_{\Lambda \cap \theta^\bot,
\nu'}$ where $\nu'$ is the restriction of $ \nu$ to
$M(\theta,\Lambda)$. Since $\ell_\nu^s$ is contained in
$\theta^\bot$ by hypothesis, the vector space $\ell_{\nu'}^s$
associated to the pair $(\nu', s)$ with respect to $\Lambda \cap
\theta^\bot$ is equal to $\ell_{\nu}^s$ and the restrictions of
$\nu$ and $\nu'$ to this subspace coincide. The inclusion (ii) holds by
the argument in the proof of (ii) $\Rightarrow$ (iii) in Proposition \ref{Prop-equiv}.
 \hfill $\Box$

\begin{Pro} \label{new}
If $1\leq k \leq d$ and $(\nu, s) \in A_k$ then the following
assertions are equivalent:
\begin{enumerate}
\item[(i)] The intersection $j_s(H^*_{\nu, \Lambda}) \, \cap \, (
\bigcup_{0 \ne \theta \leq \s} j_s (H_{\Lambda \cap \theta^\bot}))$
is empty.

\item[(ii)]  The face ${\mathcal F}_\nu$ of the polyhedron $\mathcal{N}
(\J_k) $  determined by $\nu$ is contained in the interior of  $\s^\vee$.
\end{enumerate}
\end{Pro}
\noindent{\em  Proof.}  By Proposition \ref{thetaproposition} we
have that (i) holds if and only if for any face $\theta$ of $\s$
the inclusion $\ell_\nu^s \subset \theta^\bot$ implies that
$\theta=0$, or equivalently if and only if $\ell_\nu^s \cap \int
(\s^\vee)  \ne \emptyset$.

If (ii) holds then $\ell_\nu^s \cap \int (\s^\vee)
\ne \emptyset$  by  Lemma
\ref{crucial-bis}     hence (i) holds.

Suppose that (ii) does not hold, that is there exists a vertex $w$
of  ${\mathcal F}_\nu$ which belongs to a proper face $\s^\vee
\cap \theta^\bot$ of the cone $\s^\vee$, for some $0 \ne \theta
\leq \s$. Such $w$ belongs to $\J_k$, hence it is of the form $w=
e_{j_1} + \cdots + e_{j_k}$. Since $e_{j_r}$ belongs to $\s^\vee$
it follows that $e_{j_r}$ must belong to $\theta^\bot$, for $r=1,
\dots, k$. It follows that $\ell_\nu^s \subset \theta^\bot$ by
Lemma \ref{crucial-bis}, hence (i) does not hold. \hfill $\
{\Box}$

\begin{definition} \label{Sk}
If $1 \leq k \leq d$ we define the set  $\mathcal{D}  _k $ as the
subset of cones $\t \in \bigcap_{i= 1}^k \Sigma_i$ such that the
face $\mathcal{F}_\t$ of $\mathcal{N} (\J_k)$ is contained in the
interior of $\s^\vee$.
\end{definition}

\begin{remark} \label{dd}
Notice that $\mathcal{D}_d  =\bigcap_{i= 1}^k \Sigma_i$. If $\t
\in \mathcal{D}_d$, the set $\t(d)$ is non-empty if and only if
$\stackrel{\circ}{\t}  \subset \stackrel{\circ}{\s}$.
\end{remark}

As a consequence of the results of this Section we have the
following Propositions:
\begin{Pro} \label{29}
Let us fix an integer $s_0  \geq 1$. The set $j_{s_0}
(H_\Lambda^*)\backslash\bigcup_{0 \neq \theta\leq\s}j_{s_0}
(H_{\Lambda\cap\theta^\bot})$ expresses as a finite disjoint union
of locally closed subsets, as follows:
\begin{equation} \label{trunc-s}
j_{s_0}(H_\Lambda^*)\setminus\displaystyle\bigcup_{ 0
\neq\theta\leq \s}
j_{s_0}(H_{\Lambda\cap\theta^\bot})=\displaystyle\bigsqcup_{k=1}^d
\bigsqcup_{\t\in\mathcal{D}_k} \bigsqcup_{[ (\nu, s_0) ] \in A_{k,
\t} /_{\sim}}  j_{s_0} (H_{\Lambda,\nu}^*).
\end{equation}
\hfill $\ {\Box}$
\end{Pro}

If $s_0 \geq 1$ the coefficient of $T^{s_0}$ in the auxiliary
series $P(\Lambda)$ is obtained by taking classes in the
Grothendieck ring in (\ref{trunc-s}), and then using Theorem
\ref{key-bis}.

For each cone $\t \in \mathcal{D}  _k $ we define the auxiliary series:
\begin{equation}  \label{Pktau}
P_{k,\t}(\Lambda)= (\L-1)^k \sum_{s\geq 1}   \sum_{[ (\nu, s) ]
\in A_{k, \t} /_{\sim}}   \L^{sk-\ord_{\J _k}(\nu)}T^s.
\end{equation}
\begin{Pro}
We have that
\begin{equation} \label{28}
P(\Lambda) = \sum_{k=1}^d \sum_{\t \in
\mathcal{D}_k} P_{k,\t}(\Lambda).
\end{equation}
\hfill $\ {\Box}$
\end{Pro}

\section{The proofs of the main results}

\label{racionality}

In this Section we fix a cone   $\t\in\mathcal{D}_k$ such that
$A_{k, \t} \ne \emptyset$ and we describe   the rational form of the series
$P_{k,\t}(\Lambda)$.
For convenience,  we do not stress
the dependency  on the cone $\t$ in the notations introduced in this Section.

We denote the closure of $\t(k)$ by
$\hat{\t}$. By Remark \ref{pwl}  the functions $\f_1,
\dots, \f_k$ are linear on $\t$.
More precisely,
 if $\nu_0 \in \stackrel{\circ}{\t}$  we consider the vectors $e_{i_1}, \dots, e_{i_d}$ introduced in
Proposition \ref{algo}. Then we deduce that $\f_r (\nu) = \langle \nu,
e_{i_r} \rangle $ for $1 \leq r \leq k$ and for all $\nu \in \t$, since $\t$ is a cone in
the fan $\cap_{r=1}^k \Sigma_r$.

\begin{Not}
Let us define the lattice homomorphisms
\[
\begin{array} {lcl}
\p :  \hat{N} \rightarrow \Z^{k+1} & \mbox{  by } &  (\nu, s) \mapsto (\langle \nu, e_{i_1} \rangle, \dots, \langle \nu, e_{i_k} \rangle,s),
\\
\z : \Z^{k+1} \rightarrow \Z^{2} &  \mbox{  by }  & a = (a_1, \dots, a_{k+1})  \mapsto ( k a_{k+1} - \sum_{r=1}^k a_r,  a_{k+1}) ,
\end{array}
\]
We set also  $\xi:= \z \circ \p: \hat{N} \rightarrow \Z^2$.  We
abuse of notation by denoting by the same letter the linear
extension of these maps to the corresponding real vector spaces.
\end{Not}
 Notice that the intersection of the kernel of $\pi$ (and also of
$\xi$) with the cone $\hat{\t}$ is $\{ 0 \}$.

\begin{remark} \label{rays}
If $(\nu, s) \ne (0, 0) $ belongs to a ray in $\hat{\t}$ then  $\xi ( \nu, s) \ne (0,0) $ and
\[
\xi ( \nu, s) = \left\{
\begin{array}{lcl}
(\Psi_k ( \nu ), \f_k (\nu))  & \mbox{ if } &   s = \f_k (\nu),
\\
(\Psi_{k+1} ( \nu ), \f_{k+1}  (\nu))  & \mbox{ if }  &   k \ne d  \mbox{ and } s = \f_{k+1} (\nu),
\\
(ds, s)  & \mbox{ if }  & k = d \mbox{ and }   \nu =0.
\end{array}
\right.
\]
\end{remark}
It follows from Remark \ref{rays} and Corollary \ref{c1-bis} that
 $\tilde{\t} := \xi ( \hat{\t})  \subset \R^2_{\geq 0}$ and
  $\bar{\t}:= \pi (\hat{\t}) \subset \R^{k+1}_{\geq 0} $ are strictly convex and rational for the lattices
  $\Z^2$ and $\Z^{k+1}$ respectively. Hence the map of $\C$-algebras
\[
\begin{array} {lcl}
\z_* :  \C [[ \bar{\t} \cap \Z^{k+1}    ]]  \rightarrow  \C [[ \Z^2_{\geq 0}  ]]
 = \C  [[ \L , T]] & \mbox{ given by } X^{ a} \mapsto \L^{k a_{k+1} - \sum_{r=1}^k a_r} T  ^{a_{k+1}},
\end{array}
\]
for $a = (a_1, \dots, a_{k+1})$ is well defined. If $B \subset
\bar{\t} \cap \Z^{k+1} $ the generating function $F_B := \sum_{ a
\in B} X^a$ belongs to the ring $ \C [[ \bar{\t} \cap \Z^{k+1} ]]
$.

\begin{Lem}
The sets $A_{k,\t}$,   $\bar{A}_{k, \t} := \p (A_{k,\t }) $ and
$\tilde{A}_{k,\t} : = \xi (A_{k,\t} )$ are subsemigroups, not
necessarily of finite type of $\hat{\t} \cap \hat{N}$,
$\Z^{k+1}_{\geq 0}$ and $\Z^2_{\geq 0}$ respectively. The
restriction $\p_{| A_{k,\t}}:  A_{k,\t} \rightarrow \bar{A}_{k,
\t} $ induces a bijection $ A_{k,\t} /_{\sim}  \rightarrow
\bar{A}_{k, \t}$, by $ [(\nu, s)] \mapsto \p(\nu, s)$ (see
Definition \ref{eq}). We have $ P_{k,\t}(\Lambda) = \z_*
( F_{\bar{A}_{k, \t} }). $ \label{Lemaequiv}
\end{Lem}
{\em Proof.} The result follows from Definition \ref{eq},
 Remark \ref{86} and the previous discussion. \hfill
$\ {\Box}$

\begin{Not}
 If $\r \subset \t$ is a one-dimensional cone rational for the lattice $N$ we denote by $\nu_\r$ the primitive integral vector on $\r$, that is, the generator of the semigroup $\r \cap N$.
 \end{Not}

\begin{Pro}\label{ratP_k}
If $1 \leq k \leq d- 1 $ there exists $R_{k, \t} \in \Z [\L, T]$ such that:
\begin{equation}
P_{k,\t}(\Lambda) = {R_{k, \t}}
\prod_{\r\leq\t}^{\mathrm{dim}\r=1}( 1 - \L^{\Psi_{k} (\nu_\r)}
T^{\f_{k} (\nu_\r)})^{-1}
\prod_{\r\in\Sigma_{k+1},\r\subset\t}^{\mathrm{dim}\r=1,\phi_{k+1}(\nu_\r)\neq\phi_k(\nu_\r)}(1
- \L^{\Psi_{k+1} (\nu_\r)} T^{\f_{k+1} (\nu_\r)}) ^{-1}.
\label{P-k}
\end{equation}
If $k = d$ then (\ref{P-k}) holds by replacing the term
$\prod_{\r\in\Sigma_{k+1},\r\subset\t}^{\mathrm{dim}\r=1, \phi_{k+1}(\nu_\r)\neq\phi_k(\nu_\r)}(1
- \L^{\Psi_{k+1} (\nu_\r)} T^{\f_{k+1} (\nu_\r)})$ by $ (1 - \L^d
T)$.
 Both numerator and denominator in (\ref{P-k}) are determined by the lattice $M$ and the Newton
polyhedra of the logarithmic jacobian ideals.
\end{Pro}
{\em Proof.} We call the set $\partial_-\hat\t=\{(\nu,\phi
_k(\nu))\ |\ \nu\in\t\}$ the lower boundary of $\hat{\t}$. The set
$\partial_-\hat\t $ is a convex polyhedral cone of dimension $d$.
We deduce that $A_{k, \t} = (A_{k, \t} \cap \partial_-  \hat{\t})
\sqcup   (A_{k, \t} \cap  \int(\hat{\t})) $. The sets $A_{k, \t} \cap
\partial_-  \hat{\t}$ and $  A_{k, \t} \cap  \int(\hat{\t})$
consist of the integral points for the lattice $\hat{N}$ in the
cones $\partial_-\hat\t $ and  $\hat\t $, respectively. It is easy
to see that if $(\nu, s) \in  A_{k, \t} \cap \partial_-  \hat{\t}
$ and if $(\nu', s') \in  A_{k, \t} \cap \int(\hat{\t}) $ then
$[(\nu, s) ] \ne [(\nu', s') ]$ (see Notation \ref{auxi} and
(\ref{equivalence-bis})). It follows that $\bar{A}_{k, \t} = \p
(A_{k, \t} \cap
\partial_- \hat{\t}) \sqcup \p ( A_{k, \t} \cap  \int(\hat{\t}))
$. We set $\bar{A}_{k, \t}^\circ: = \p ( A_{k, \t} \cap
\int(\hat{\t})) $ and $\bar{A}_{k, \t}^- := \p (A_{k, \t} \cap
\partial_- \hat{\t}) $. It follows that
\begin{equation} \label{ktau}
P_{k,\t}(\Lambda) = (\L -1) ^k  \,(  \z_* ( F_{\bar{A}_{k,
\t}^\circ
 } ) +   \z_* ( F_{\bar{A}_{k, \t}^- })) .
\end{equation}

The semigroups ${\bar{A}_{k, \t}^\circ }$ and ${\bar{A}_{k, \t}^-
}$ are the images by $\p$ of the semigroups of integral points in
the relative interiors of the cones $\hat{\t}$ and
$\partial_-\hat{\t}$, respectively. We apply Theorem
\ref{teoremaapendice} (see Section \ref{app}) using that the
kernel of $\p$ intersects the cone $\hat{\t}$ only at $0$.

It follows that the denominator of the rational form of
$F_{\bar{A}_{k, \t}^\circ
 }$  (resp.  of $F_{\bar{A}_{k, \t}^- }$)
is the product of terms $1 - X^{\p (b) } $, for $b$ running
through the primitive integral vectors in the edges of
 $\hat{\t}$ (resp. of $\partial_-\hat\t $)
 while the numerator is a polynomial in $\Z [\bar{\t} \cap \Z^{k+1}]$.
The rational form of $P_{k,\t}(\Lambda)$ is the image of the rational form of $F_{\bar{A}_{k, \t}^\circ
 }  +   F_{\bar{A}_{k, \t}^- }$  by the  homomorphism $\z_*$
since     the image by $\z_*$ of the denominator does not vanish by    Remark \ref{rays}.

 If $1 \leq k \leq d-1$ and if $\r$ is an edge of $\hat{\t}$ which is not contained in $\partial_- \hat{\t}$
then it is necessarily of the form $\r = ( \nu, \f_{k+1} (\nu))$ for $\nu \in \stackrel{\circ}{\t}$ in some edge of $\Sigma_{k+1}$. If
$k = d$ the only edge of $\hat{\t}$ which is not  contained in $\partial_- \hat{\t}$ is $(0,1) \R_{\geq 0}$.
 Finally by this discussion and Remark
\ref{rays} the denominator of this rational form is as indicated in
(\ref{P-k}).
 \hfill $\Box$

\begin{remarks} $\,$ \label{corPolo1-LdT}
\begin{enumerate}
\item[(i)] For $k=1,\ldots,d-1$ and $\t\in\mathcal D_k$, the
factor $\ 1-\L^dT$ does not appear in the denominator of
$P_{k,\t}(\Lambda)$.

\item[(ii)]  The factors in the denominator of the rational form
(\ref{P-k}) of $P_{k,\t}(\Lambda)$ are of the form $1 - \L^a T^b$
with $(a, b) \in B(\Lambda)$. We use that $\cup_{i=1}^k
\Sigma_i^{(1)}$ is the set of rays in the fan  $\cap_{i=1}^k \Sigma_i$
and Definition \ref{Sk}.

\item[(iii)] The term $1-\L^dT$ appears in the denominator of $P_{d, \t}
(\Lambda)$ with multiplicity one.
\end{enumerate}
\end{remarks}

{\em Proof of Proposition \ref{P_curvasmonomiales}.} If $d= 1$
then the toric variety $Z^{\Lambda}$ is a monomial curve.
 Let $\nu_0$ be the generator
of the semigroup $\s \cap N \cong \Z_{\geq 0}$. The monomial curve
$Z^\Lambda$ is parametrized by $x_i=t^{m_i} $, where $m_i :=
\langle\nu_0,e_i\rangle $, for $i=1,\ldots,n$. The multiplicity of
$Z^\Lambda$ at $0$ is $m=\min_{i =1, \dots , n
}\{\langle\nu_0,e_i\rangle\} = \mbox{\rm ord}_{\J_1} (\nu_0)$.

By Definition \ref{ak},  the set $A_1$ is  $A_1 = \{ (\nu, s) \mid  \nu =
r \nu_0, \, \,  0 < m r \leq s  \}$. By Theorem \ref{key-bis} it
follows that $P(\Lambda) = \sum_{0 < r m \leq s} (\L -1) \L^{s- r
m} T^s = \frac{\L-1}{1-\L
    T}\frac{T^m}{1-T^m}$.
\hfill\ $\Box$

\medskip

 {\em Proof of Theorem  \ref{PLambdaRac}.}
The lattice $M$ and the Newton polyhedra of the ideals
     $\mathcal J_k$  determine and are determined by duality by
$N$ and the functions  $\mbox{\rm ord}_{\mathcal J_k  } $, for   $k=1, \dots, d$.
The proof  follows from
Propositions   \ref{ratP_k}, Formula
 (\ref{28}) and Remark \ref{corPolo1-LdT}.
      \hfill$\Box$

\medskip

{\em Proof of Corollary \ref{P-geom}.} It is a consequence of
Theorem \ref{PLambdaRac}, Proposition \ref{descompPgeom} and
 Example \ref{P-sigma}.  \hfill $\Box$

\medskip

{\em Proof of Corollary \ref{P-locales}.} Nicaise observed in
\cite{Nicaise} that the motivic Poincar\'e series of an affine
normal toric variety $Z^\Lambda$ has an expansion in terms of the
local motivic Poincar\'e series at the distinguished points of the
orbits, namely:
\[
P_{\mathrm {geom}}^{Z^\Lambda}(T)=\displaystyle\sum_{\theta
\leq\s}(\L-1)^{\mbox{\rm\tiny codim}\theta}P_{\mathrm {geom}}^{(Z^\Lambda,o_\theta)}(T).
\]
For each $\theta \leq \s$ there exists an open set of $Z^{\Lambda}$
containing the distinguished point $o_\theta$  of the orbit
$\O_{\theta}^{\Lambda}$, which is isomorphic to
$\O_{\theta}^{\Lambda} \times Z^{\Lambda_\theta}$, where
$\Lambda_\theta= \s^\vee_\theta \cap M_\theta$ is the image of $\Lambda$ by the canonical map $M
\rightarrow M / M \cap \theta^\bot$  (see \cite{Fu} page 29). It
follows that the germ $(Z^\Lambda, o_\theta)$ is analytically
isomorphic to $(Z^{\Lambda(\theta)}, 0) $.
 \hfill$\Box$

\medskip
 {\em Proof of Corollary \ref{poles-surfaces}.}  It follows from Theorem \ref{PLambdaRac}
 and Example \ref{P-codim1} by taking into consideration that if $\r $ is a ray of $\Sigma_1$
and if $\nu_\r \in \stackrel{\circ}{\s} $ then we get $\f_1
(\nu_\r) = \f_2 (\nu_\r)$, thus $\Psi_1 (\nu_\r) = \Psi_2 (\nu_\r)
=0$.  \hfill $\Box$

{\em Proof of Proposition \ref{val3}.} It follows from the
definitions by using  (\ref{menor-bis}) and
(\ref{val2}).  \hfill $\Box$

\section{Motivic volume of a toric variety} \label{volume}

 We give a formula  for the motivic volume    $\mu (H_\Lambda)$
of the space of arcs $H_\Lambda$ of the  toric variety $Z^\Lambda$
 in terms of
the support function $\ord_{\J_d}$.
This
formula generalizes the one given in \cite{LR}.

      If $\t \in \Sigma_d$ we
denote by $\eta_\t:\Z[\t\cap N]\longrightarrow\Z[\L^{\pm 1}]$ the
ring homomorphism defined by $\eta_\t
(x^\nu)=\L^{-\ord_{\J_d}(\nu)}$. The generating function $
F_{\stackrel\circ\t\cap N} $ has a rational form
$ F_{\stackrel\circ\t\cap N} = R_{\stackrel\circ\t\cap N} \prod_{\r \leq \t}^{\dim \r =1} ( 1 - x^{\nu_\r})^{-1}$,
for some    $ R_{\stackrel\circ\t\cap N} \in \Z [\t \cap N]$
(see Proposition
\ref{open-cone}).

\begin{Pro}
\[ \mu(H_\Lambda)= (\L -1)^d \sum^{\t\in\Sigma_d}_{{\stackrel\circ\t\cap\stackrel\circ\s\neq\emptyset}}
\eta_\t (R_{\stackrel\circ\t\cap N})
\prod_{\r \leq \t}^{\dim \r =1} ( 1 - \L^{- \ord_{\J_d} ( \nu_\r) })^{-1}
\] \label{vol-mot2}
\end{Pro}
{\em Proof.} By Theorem \ref{ThVolMot} (see \cite{DL1}) we have
that the limit
$\mu(H_\Lambda)=\lim_{m\to\infty}[j_m(H_\Lambda)]\L^{-md}$
converges in $\hat{\mathcal M}$.  We deduce that
$\mu(H_\Lambda)=((1-\L^d T)P_{\mathrm {geom}}^{(Z^\Lambda,0)}(T))_{|T= \L^{-d}}$,
by comparing with the
definition
of the series $P^{(Z^\Lambda,0)}_{\mathrm {geom}}$, taking into
account that $1-\L^dT$ is a simple pole of $P_{\mathrm
{geom}}^{(Z^\Lambda,0)}$ by Propositions \ref{descompPgeom} and
Remark \ref{corPolo1-LdT}.
By Proposition \ref{ProVolMot}  the equality
$\mu(H_\Lambda)=\mu(H_\Lambda^*)$ holds. By Remark \ref{corPolo1-LdT}
the term $1-\L^dT$ does not define a pole of $P_k(\Lambda)$ for
$k=1,\ldots,d-1,$ hence
$\mu(H_\Lambda)=((1-\L^dT)P_d(\Lambda))_{|_{T=\L^{-d}}}$ and
\[
\mu(H_\Lambda)=\sum^{\t\in\Sigma_d}_{{\stackrel\circ\t\cap\stackrel\circ\s\neq\emptyset}}
\sum_{\nu\in \stackrel\circ\t \cap N}(\L-1)^d\L^{-\ord_{\J_d}(\nu)}
= (\L -1)^d \sum^{\t\in\Sigma_d}_{{\stackrel\circ\t\cap\stackrel\circ\s\neq\emptyset}}
\eta_\t (F_{\stackrel\circ\t\cap N}(x)).\]
  Notice that $\eta_\t (x^{\nu_\r})\neq 1$ if $\nu_\r$ is
a primitive vector in a ray of $\r \in\Sigma_d$ (see Remark
\ref{dd}).     We deduce that
  $ \eta_\t (F_{\stackrel\circ\t\cap N}) =  \eta_\t (R_{\stackrel\circ\t\cap N})
\prod_{\r \leq \t}^{\dim \r =1} ( 1 - \L^{- \ord_{\J_d} ( \nu_\r) })^{-1}$.
\hfill $ \Box$

We deduce a formula for the motivic volume of the space of arcs
$H(Z^\Lambda)$ of $Z^\Lambda$ (without fixing the origin of the
arcs) in terms of the local data, as a consequence of Proposition
\ref{vol-mot2} and Corollary \ref{P-locales}.  The same formula
also holds if $Z^{\Lambda}$ is locally analytically unibranched
(see \cite{tesis}).

\begin{Pro} \label{vol-motGlobal2}
If $Z^\Lambda$ is  an affine normal toric variety then we have
that
\[
\mu(H(Z^\Lambda))=\sum_{\theta\leq\s}  \, (\L-1)^{\mbox{{\rm \small codim}}\theta} \, \mu(H_{\Lambda(\theta)}).
\]            \hfill $ \Box$
\end{Pro}

\section{Geometrical definition of the logarithmic jacobian
ideals} \label{tor}

 We introduce the geometrical definition of the $k^{th}$-{\em logarithmic jacobian ideal} of an affine toric variety
$Z^\Lambda$ of dimension $d$ for $1\leq k\leq d$, following
\cite{Oda} Chapter 3, and \cite{LR} Appendix. We denote by $D$ the
equivariant Weil divisor defined by the orbit closures of
codimension one in $Z^{\Lambda}$. We denote by
$\Omega_\Lambda$ the $\C[\Lambda]$-module of K\"ahler differential
forms of $Z^\Lambda$ (over $\C$). The module $\Omega_{\Lambda}
(\log D)$  of one forms on $Z^{\Lambda}$ with logarithmic poles
along $D$ is identified with $\C[\Lambda ] \otimes_\Z M$ and we
have a map of $\C[ \Lambda]$ modules:
\[
 \varphi:  \Omega_{\Lambda} \rightarrow    \Omega_{\Lambda} (\log D), \quad dX^\g \mapsto X^\g \otimes \g, \mbox{ for }
\g \in        \Lambda.
\]

If $1 \leq k \leq d$ we set
$
 \y^k \varphi :      \Omega_{\Lambda}^k    \rightarrow   \Omega_{\Lambda}^k (\log D), \quad dX^{\g_1} \y \cdots \y   dX^{\g_k} \mapsto X^{\g_1 + \cdots + \g_k} \otimes
(\g_1 + \dots + \g_k)$,
for $\g_i $ in $\Lambda$, where $\Omega_{\Lambda}^k (\log D)$ is
identified with $\C[\Lambda] \otimes_\Z \bigwedge^k M$. For $k =d$,
fixing a basis $u_1, \dots, u_d$ of the lattice $M$
provides an isomorphism
\[
  \f:  \wedge^d M \rightarrow \Z, \quad u_1 \y \dots \y u_d \mapsto 1,
\]
which is, up to sign, independent of the choice of the basis.

\begin{definition} The
         {\em $k^{th}$-logarithmic jacobian ideal}  of $Z^\Lambda$
is the ideal of $\C[\Lambda]$ generated  by  the set  $\f
(\y^k\varphi ( \Omega_\Lambda^k)\y\bigwedge^{d-k}M)$.
\end{definition}

\begin{Pro}
The $k^{th}$-logarithmic jacobian ideal of $Z^\Lambda$ is the monomial
ideal $\J_k$ defined by (\ref{j-k}), for $k=1, \dots, d$.
\end{Pro}
{\em Proof.} The proof follows from the definitions since
$\Omega_\Lambda$ is generated by $dX^{e_i}$, for $e_1, \dots, e_n$
generators of the semigroup $\Lambda$. \hfill $\ {\Box}$

                                       The {\em  Nash blowing up} $v: N_S \rightarrow S$
of an algebraic variety    $S$
is the minimal proper birational map such that $v^* \Omega^1_S$  has
a locally free quotient of rank $\dim S$.   The fibers of $v$ at a point $x \in S$ are
the limiting positions of tangent spaces at smooth points of $S$ tending to the point $x$.

 \begin{Pro}     {\rm (\cite{GN,LR,Teissier})}
The blowing up of $Z^{\Lambda}$ with center its $d^{th}$-logarithmic
jacobian ideal $\J_d$  is the Nash blowing up of
$Z^{\Lambda}$.     \label{NashZ}
\end{Pro}

\section{Generating functions of projections of subset of cones}   \label{app}

In this Section
we state some auxiliary results on the generating function of certain subsets of integral points in a rational polyhedral cone.
See \cite{Brion, BP, Stanley} for an expository papers on this and related subjects.
The content of this section is independent of the rest of the paper.

Let $N \subset \R^d$ be a rank $d$ lattice and $\t$ strictly convex cone
rational for the lattice $N$.

\begin{definition}
The {\em generating
function} of a set $B \subset \t \cap N$ is the series
$F_B(x)=\sum_{a\in B}x^a \in \Z [[\t \cap N]]$.   The series   $F_B(x)$ is {\em rational}  if
there exist $p(x), q (x) \in \Z[ \t \cap N ]$ such that $q(x) F_B (x) =
p(x)$. In that case the ratio $p(x)/q(x)$ is well-defined and it
is called the {\em sum} of the series $F_B(x)$.
\end{definition}

  We denote by $\nu_\r$ the generator of  the semigroup $\r
\cap N$ for each edge $\r$ of $\t$.
   The following Proposition is well-known    (see \cite{Stanley} Section 4.6).
\begin{Pro}  \label{open-cone}
The generating function
$F_{\stackrel{\circ}{\t}\cap N}(x) $  is of the form:
\begin{equation}   \label{funcGeneratriz}
    F_{\stackrel{\circ}{\t}\cap N}(x)= R_{\stackrel{\circ}{\t}\cap N}
\displaystyle\prod_{\r \leq \t, \dim \r = 1} (1-x^{\nu_\r})^{-1},
\mbox{ with }  R_{\stackrel{\circ}{\t}\cap N}   \in \Z[\t \cap N].
\end{equation}
\end{Pro}
\begin{remark}      \label{rem-rho}
The statement of Proposition \ref{open-cone} remains true
 if we replace the vector $\nu_\r$ by a non-zero vector in $\r \cap N$
for each edge $\r$ of $\t$.
\end{remark}

   Let  $\p : N \rightarrow \Z^r$ be a map of lattices for some $1 \leq r \leq d$.
We abuse of notation by denoting with the same letter the extension of $\p$ to a map of real vector spaces $  N_\R \rightarrow  \R^r$.
We suppose that
  \begin{equation}
     \t \cap \mbox{ ker } \p    =\{ 0\}.
    \label{ele2}
  \end{equation}
This condition implies that the cone $\bar{\t}: = \p (\t)$ is strictly convex.
For simplicity we set $A: = \stackrel\circ\t\cap N$ and $\bar{A}:= \p (A)$.
The sets $A$ and $\bar A$ are subsemigroups, not necessarily of finite type,
 of $\t \cap N$ and $\bar{\t} \cap \Z^r$ respectively.
Notice that for each edge $\bar \r$ of $\bar{\t}$ there exists at least one
edge ${\r}$ of $\t$ such that $\p (\r) = \bar \r$, hence we have that
\begin{equation}         \label{rho2}
 0 \ne \p (\nu_\r)  \in \bar{\r} \cap \Z^r.
\end{equation}

\begin{The} \label{teoremaapendice}
The
generating function $F_{\bar{A}} (x) $ of
$\bar{A}   $ is of the form:
\begin{equation}           \label{FA}
       F_{\bar{A}} (x) = R_{\bar{A}} \prod_{{\r} \leq {\t}, \,
\dim {\r} =1} (1-x^{\p ({\nu}_{{\r}}) })^{-1}, \mbox{ with
some } R_{\bar{A}} \in\Z[\bar{\t} \cap \Z^r].
\end{equation}
\end{The}

We introduce some notations and results before proving Theorem \ref{teoremaapendice}.

If $a \ne b \in N_\Q$  we define the {\em length with respect to
$N$} of the segment  joining $a$ and $b$ by
 $\textrm{lg}(a,b) : = r $ if  $a-b=rc$ where $r \in \Q_{>0}$ and $c \in N$ is
 a primitive vector.

\begin{figure}[htbp]
$$\epsfig{file=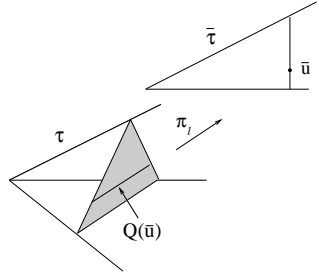, height=36mm}$$
\caption{The shaded region is the preimage by $\p_\ell$ of the
segment containing $\bar u$ in the Figure.\label{apendice2}}
\end{figure}

We denote by $\mathcal A$ and $\mathcal B$ the following sets:
\[\mathcal A=\{\bar\r\mbox{ edge of }\bar\t\ |\ \dim\pi_\ell^{-1}
(\bar\r)\cap\t=1\}\ \mbox{ and }\ \mathcal B=\{\bar\r\mbox{ edge of }
\bar\t\ |\ \dim\pi_\ell^{-1}(\bar\r)\cap\t>1\}.\]
\begin{Lem}
If $r < d$ and  $\bar\r\in\mathcal B$ then there exists $\bar
u_0\in\bar\r\cap\Z^r$ such that $(\bar u_0+ {\rm{int}}(\bar\t))\cap\bar
N\subset\bar{A}$.
\label{PropB}
\end{Lem}
{\em Proof.}
 If $\bar u\in\bar\t$ we denote by $Q(\bar u)$ the set
$Q(\bar u):=\pi_\ell^{-1}(\bar u)\cap\t$
(see Figure  \ref{apendice2}).
If $\bar u \in \bar\t \cap \Z^r $ then $Q(\bar u)$ is a
 rational polytope  for the lattice $N$.
We denote by $\r(\bar u)$ the ray spanned by the sum of all the vertices of
the polytope $Q(\bar u)$,  by $ b(\bar u)$ the vector
$Q(\bar u)\cap\r(\bar u)$ and by $\d(\bar u)$ the number
$\d(\bar u) := \max\{\mathrm{lg}(v,{b}(\bar u)) \}$, for $
v$ running through the vertices of  $Q(\bar u) $.
Notice that if $t \in \R _{\geq 0} $ then we have that $\r(t\bar u)=
    \r(\bar u)  $ and
        \begin{equation}
    Q(t \bar u)=tQ(\bar u) ,   \quad   {b}(t\bar u)=t{b}(\bar u),  \quad   \d (t\bar u)=t \d (\bar u).
        \label{lin-B}
        \end{equation}
 If $\bar\r \in \mathcal{B}$ and  $\bar u
\in\bar\r\cap\Z^r$ then the polytope $Q(\bar u)$ is of dimension $\geq 1$ by definition of $\mathcal B$,
hence   $\d(\bar u) >0$.
Let $\bar u_0$ be a vector in $\bar\r\cap\Z^r$ such that $\delta(\bar u_0)>1$.

A vector $\bar{w}  \in (\bar{u}_0 +  {\rm{int}}(\bar{\t}))\cap\Z^r  $ is of the form
   $\bar{w} := \bar{u}_0 + \bar{v}$, with   $\bar{v} \in  {\rm{int}}(\bar{\t})\cap\Z^r$.
 By linearity of $\p$ we have the inclusion
\begin{equation}
Q(\bar u_0) + Q (\bar v ) \subset Q ( \bar w)   .
\label{nov}
\end{equation}
Since $\p$ maps $ {\rm{int}}(\t)$ onto $ {\rm{int}}(\bar\t)$ the polytopes $Q(\bar v)$
and $Q(\bar w)$ are of dimension $d-r \geq 1$. Thus, there exists a vector $v \in N_\Q$  such that
$v  \in{ {\rm{int}}}(Q(\bar v))$,
where int denotes the relative interior. By (\ref{nov}) we obtain
$v+{ {\rm{int}}}(Q(\bar u_0))\subset { {\rm{int}}}(Q(\bar w))$.
Since $\delta(\bar u_0)>1$ it follows that $Q(\bar w)$ contains
 a rational segment $v+I$ of integral length $>1$. Hence $Q(\bar w)$
contains a point $w$ of the lattice $N$ in the set $ {\rm{int}}(Q(\bar w))\subset {\rm{int}}(\t)$.
It follows that $\bar w=\p(w)$ belongs to $\bar{A}$. \hfill $\Box$

\medskip

{\em Proof of the Theorem \ref{teoremaapendice}}. We deal first
with  the case of a simplicial cone $\bar\t$.   In this case we denote
$\mathcal A=\{\bar\r_1,...,\bar\r_a\}$ and $\mathcal
B=\{\bar\r_{a+1},...,\bar\r_{a+b}\}$, where $a, b \geq 0$ and $a+
b = \dim\bar \t$.

If $\bar \r_j \in \mathcal B$ we denote by
$\bar u_j$ the vector $\bar u_j \in \bar \r_j \cap \Z^r $ such
that $(\bar u_j +  {\rm{int}}(\bar \t)) \cap\Z^r \subset \bar A$ (see
Lemma \ref{PropB}). The set $S: =\bigcup_{\r_j \in \mathcal B} (\bar u_j+ {\rm{int}}(\bar\t))\cap\Z^r$ is
contained in $\bar{A}$. If
$S':=\bar{A}\setminus S$ then $F_{\bar{A}} (x) = F_{S} (x)
+ F_{S'} (x)$.

We deal first with the rational form of $F_{S} (x)$.
If  $\emptyset \ne J \subset \mathcal B$
we set $R_J := \bigcap_{\r_j \in J} (\bar
u_{j}+ {\rm{int}}(\bar\t))\cap\Z^r$. Notice that $R_J = (\bar{u}_J + {\rm{int}}(\bar{\t}))\cap \Z^r$,
where $\bar u_J := \sum_{\r_j \in J} \bar u_j$.  By Proposition \ref{open-cone}
the series
$F_{R_J}(x)=x^{\bar u_J}F_{ {\rm{int}}(\bar\t)\cap\Z^r}(x)$  is of rational form.
By Remark \ref{rem-rho}  and Formula (\ref{rho2})
its denominator can be taken as in (\ref{FA}).
By the
inclusion-exclusion principle we deduce that $F_S (x)$
has a rational form as indicated in the statement of
Theorem         \ref{teoremaapendice}.

                 If $\bar \r_j \in \mathcal A$ we denote by $\r_j$ the edge of $\t$ such that
$\p (\r_j) = \bar \r_j$, for $j=1, \dots, a$.
To study the rational form of $F_{S'}(x)$ we set
\[G=\{\sum_{i=1}^a \l_i \p (\nu_{\r_i}) +\sum_{j=1}^b\mu_j \bar u_{a+j}\ |\ 0<\l_i,\mu_j\leq 1,\mbox{ for }
 1\leq i\leq a,\ 1\leq j\leq b\}.\]
      If $\vec n = (n_1, \dots, n_a) \in\Z_{\geq 0}^a$ we denote by $C_{\vec n}$ the set
$C_{\vec n}:=   n_1 \p (\nu_{\r_1})     +\cdots+n_a\p (\nu_{\r_a}) +G$
and by $k_{\vec n}$ the integer $k_{\vec n}: =\# (C_{\vec n}\cap\bar{A})$, where $\#$ denotes cardinal.
Then we have a partition
\begin{equation}     \label{partition}
         S' =
\bigsqcup_{\vec{n} \in \Z^a_{\geq 0}}    C_{\vec n} \cap \bar{A} .
\end{equation}

If $\vec n\in\Z_{\geq 0}^a$ we have the bound:
\begin{equation}
k_{\vec n}\leq \# (G\cap\Z^r).
\label{ecCard}
\end{equation}
Denote by $\{\vec e_1,\ldots,\vec e_a\}$ the canonical basis of $\Z^a$. We have that
\begin{equation}
k_{\vec n}\leq k_{\vec n+\vec e_j}\ \ \mbox{ for }1\leq j\leq a\mbox{ and }\vec n\in\Z_{\geq 0}^a.
\label{ecSucc}
\end{equation}
We deduce from (\ref{ecCard}) and (\ref{ecSucc}) that there exists $m\in\Z_{>0}$ such that
$
k_{\vec n}=k_0$
for all $\vec n \in \Z^a_{\geq 0}$
such that      $n_j \geq m$ for some   $1 \leq j \leq a$.
If $\vec n$ verifies this condition and if $1 \leq i \leq a$  we have  the equality:
\begin{equation}
 C_{\vec n}\cap\bar{A}+
\p ( \nu_{\r_i} )
=(C_{\vec n+\vec e_i})\cap\bar{A}\mbox{ for }1\leq i\leq a.
 \label{ecPeriod}
 \end{equation}
 We deduce from these observations and (\ref{partition}) that \[
        F_{S'}(x)=\sum_{ \forall 1 \leq i \leq a  : 0 \leq n_i < m}   \, \,
\sum_{\bar\nu\in C_{\vec n}\cap\bar{A}}x^{\bar\nu}+ \sum_{ \exists
1 \leq i \leq a \,  : m \leq n_i  }       \, \, \sum_{\bar\nu\in
C_{\vec n}\cap\bar{A}}x^{\bar\nu}.
                          \]
The first term is a finite sum, while the second is of the form $R(x)\prod_{j=1}^a(1-x^{\p (\nu_{\r_j} )})^{-1}$
for some $R(x)\in\Z[x]$, by (\ref{ecPeriod}) and a similar argument as the one used for $F_S(x)$.

If $\bar\t$ is not simplicial, let $\bar\Sigma$ be a simplicial subdivision of $\bar\t$
such that every edge of $\bar\Sigma$ is an edge of $\bar\t$ (see Chapter V, Theorem 4.2, page 158 \cite{Ewald}).
If $\bar\theta\in\bar\Sigma$ the set $\theta:=\p^{-1}( \bar{\theta} )\cap\t$ is a rational cone
for the lattice $N$ and  $\theta\cap \mbox{{\rm ker}} (\p) =(0)$. We have that
$\p( {\rm{int}}(\theta))= {\rm{int}}(\bar\theta)$ and  $\p( {\rm{int}}(\theta)\cap N)=\bar{A}\cap {\rm{int}}(\bar\theta)$.

By the assertion in the simplicial case $F_{\bar{A}\cap {\rm{int}}(\bar\theta)}(x)$
 has a rational form as in the statement of the Theorem. The result follows since
$F_{\bar{A}}(x)=\sum_{\bar\theta\in\bar\Sigma}F_{\bar{A}\cap {\rm{int}}(\bar\theta)}(x)$.
\hfill$\Box$

\section{An Example}
\label{toric-example}

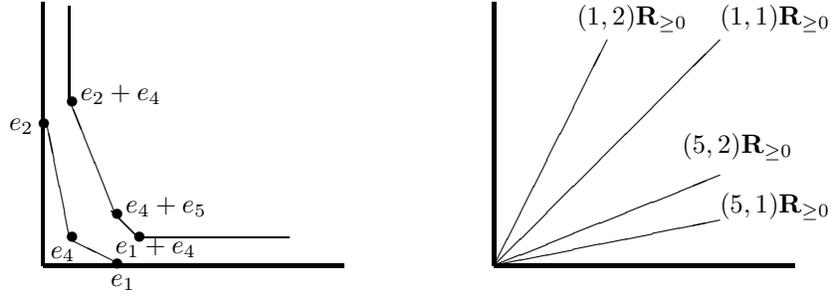
\begin{figure}[h]
\unitlength=1mm
\begin{center}
\begin{picture}(100,35)(0,0)
\linethickness{0.5mm}

\put(0,0){\line(0,1){35}} \put(0,0){\line(1,0){40}}
\put(60,0){\line(0,1){35}} \put(60,0){\line(1,0){40}}
\put(9,-0.7){$\bullet$}\put(9,-2.7){$e_1$}
\put(-0.7,18){$\bullet$}\put(-4.5,18){$e_2$}
\put(3,3){$\bullet$}\put(1,1){$e_4$}
\linethickness{0.15mm}
\put(9.5,0.5){\line(-2,1){6}}
\put(3.6,3.5){\line(-1,5){3}}

\put(12,3){$\bullet$} \put(9.6,1.6){$e_1+e_4$}
\put(9,6){$\bullet$} \put(11,7){$e_4+e_5$}
\put(3,21){$\bullet$} \put(5,22){$e_2+e_4$}

\put(12.7,3.7){\line(1,0){20}}
\put(12.7,3.7){\line(-1,1){3.6}}
\put(9.7,6.7){\line(-2,5){6}}
\put(3.5,21.5){\line(0,1){13}}

\put(60,0){\line(1,1){30}}\put(90,32){$(1,1)\R_{\geq 0}$}

\put(60,0){\line(1,2){15}}\put(71,32){$(1,2)\R_{\geq 0}$}

\put(60,0){\line(5,2){30}}\put(85,15){$(5,2)\R_{\geq 0}$}

\put(60,0){\line(5,1){30}}\put(90,7){$(5,1)\R_{\geq 0}$}

\end{picture}
\caption{The Newton polygons of $\J_1$ and $\J_2$ and the subdivision $\Sigma_1\cap\Sigma_2$     \label{h}}
\label{lala}

\end{center}
\end{figure}

We consider the semigroup $\Lambda$ generated by $
e_1=(3,0),e_2=(0,6),e_3=(5,0),e_4=(1,1),e_5=(2,1)$ and $
e_6=(1,4)$. With notations of Section \ref{sec-tor} the  cone $\s$
is $\R^2_{\geq 0}$, the lattice $M$ is $\Z^2$ and the semigroups
$\Lambda\cap\theta_i^\bot$ for $i=1,2$ are
$\Lambda\cap\theta_1^\bot=(0,6)\Z_{\geq 0}$ and
$\Lambda\cap\theta_2^\bot=(3,0)\Z_{\geq 0}+(5,0)\Z_{\geq 0}$.
The Newton polyhedron of $\J_1$ (resp. of $\J_2$) has vertices $e_1$, $e_4$ and $e_2$
(resp. $e_1 + e_4$, $e_4 + e_5$ and $e_2 + e_4$), see Figure \ref{h}.

The surface $Z^{\Lambda}$ is defined by binomial equations in $\A_\C^6$.
The normalization of the germ $(Z^{\Lambda},0)$ is smooth.
Notice that in this case there is no Hirzebruch-Jung
data from the minimal resolution
(compare with the results of \cite{LR}  in the case of a normal toric surface singularity).

By  Proposition \ref{descompPgeom}  $P_{\mathrm {geom}}^{(Z^\Lambda,0)}(T)$ is equal to
\begin{center}
$P_{\mathrm {geom}}^{(Z^\Lambda,0)}(T) = (1-T)^{-1} +
P(\Lambda\cap\theta_1^\bot) + P(\Lambda\cap\theta_2^\bot)  + P(\Lambda)$.
\end{center}
   By   Proposition \ref{P_curvasmonomiales}
we get
$P(\Lambda\cap\theta_1^\bot)=\frac{\L-1}{1-\L T}\frac{T}{1-T}$ and
$P(\Lambda\cap\theta_2^\bot)=\frac{\L-1}{1-\L
T}\frac{T^3}{1-T^3}$.
The sum of $P(\Lambda)$ is of the form  $Q_\Lambda \prod_{(a, b) \in B(\Lambda)} (1 - \L^a T^b)^{-1}$ where
$Q_\Lambda  \in \Z[\L , T]$. The set
 $B(\Lambda) = \{ (2,1)$, $(1,1)$, $(0,3)$,  $(0,6)$, $(1,3)$, $(5,12) \}$
is determined easily from the table below in which
 we give the values of the functions $\phi_1$, $\phi_2$
and $\Psi_2$ for the primitive vectors in the rays $\r$ of
$\Sigma_1  \cap \Sigma_2 $.
\vspace{2mm}

\centerline{\begin{tabular}{|c|c|c|c|c|c|c|}
  \hline
    & $(1,0)$ & $(0,1)$ & $(1,2)$ & $(5,1)$ & $(1,1)$ & $(5,2)$\\
  \hline
  $\phi_1$ & 0 & 0 & 3 & 6 & 2 & 7 \\
  \hline
  $\phi_2$ & 1 & 1 &  3 & 6 & 3 & 12 \\
  \hline
  $\Psi_2$ & 1 & 1 & 0 & 0 & 1 & 5 \\
  \hline
\end{tabular}}

\vspace{2mm}

  We have that  $\mathcal{D}_1 = \{ \t \}$ where $\t = (1,2)\R_{\geq 0}+(5,1)\R_{\geq 0} \in
\Sigma_1$  (see Definition \ref{Sk}). We have that $\mathcal{D}_2 =   \Sigma_1 \cap \Sigma_2$ and
$
 P(\Lambda) = P_{1, \t} (\Lambda) +  P_2(\Lambda)$, where
$P_2(\Lambda) :=  \sum_{\theta \in \Sigma_1 \cap \Sigma_2,
\stackrel{\circ}{\theta} \subset \stackrel{\circ}{\s}}  P_{2, \theta}   $.
We determine     the rational form  of  $ P_2(\Lambda)$
by computing first  the rational form of the
generating series $F_{\stackrel{\circ}{\theta} \cap N }$,  for $\theta \in \Sigma_1 \cap \Sigma_2$ with
$\stackrel{\circ}{\theta} \cap \stackrel{\circ}{\s} \ne \emptyset$. Then we apply
to each term  $F_{\stackrel{\circ}{\theta} \cap N}$ a suitable monomial transformation
(see the proof of  Propositions \ref{ratP_k} and \ref{open-cone}).
We check that
\begin{center}
$P_2(\Lambda)=\frac{(\L-1)^2}{1-\L^2T}\Big( \frac{\L T^4}{(1-\L
T)(1-T^3)}+ \frac{T^3}{1-T^3}+\frac{\L T^6}{(1-T^3)(1-\L
T^3)}+\frac{\L T^3}{1-\L T^3} + $
$+\frac{\L^2T^5+\L^4T^{10}+\L^6T^{15}}{(1-\L
T^3)(1-\L^5T^{12})}+\frac{\L^5T^{12}}{1-\L^5T^{12}}+\frac{\L^2T^6+
\L T^6+\L^4T^{12}+\L^3T^{12}+\L^5T^{18}}{(1-\L^5T^{12})(1-T^6)}+$
$\frac{T^6}{1-T^6}+\frac{\L T^7}{(1-T^6)(1-\L T)}\Big)$.
\end{center}

\begin{figure}[htbp]
$$\epsfig{file=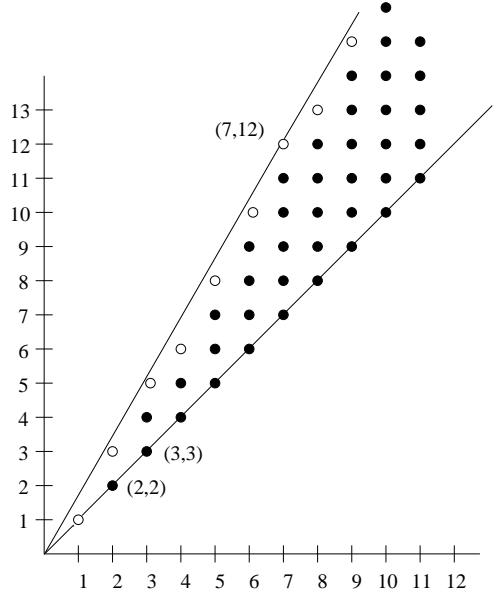, height=78mm}$$
\caption{The black  (resp. white) circles denote  elements of $\bar{A}$
(resp. of  $(\bar\t  \cap \Z^2 ) \setminus  \bar{A}$).} \label{semigrupo}
\end{figure}

We   determine the term $P_{1,
\t} (\Lambda)$ (see  Proposition \ref{ratP_k}).  The cone  $\t$ is associated
to the vertex $e_4$ of $\mathcal{N} (\J_1)$ and it is  subdivided
by $\Sigma_2$ with the rays $\r_1=(5,2)\R_{\geq 0}$ and
$\r_2=(1,1)\R_{\geq 0}$.
We describe first the
generating function $F_{\bar{A}} (x)$ of the semigroup  $\bar{A}=\{(\phi_1(\nu),s)\in\Z^2\ |\ \nu\in
\stackrel{\circ}{\t} \cap N,\ \phi_1(\nu)\leq s<\phi_2(\nu)\}$ (see  Figure \ref{semigrupo}).

We have that $\bar{A}$ is a subsemigroup of $\bar{\t} \cap \Z^2$ where
$\bar{\t} = \R_{\geq 0} (1,1) + \R_{\geq 0} (7,12)$.        We set
$G' := \{ (0,0), (2,3), (3,5) , (5,8), (6,10) \}$  and $G : = G' \cup \{  (1,1), (4,6) \}$.
We have the partitions:
\begin{center}
     $(\bar{\t} \cap \Z^2 )\setminus \bar{A}  =  \bigsqcup_{p \geq 0}  G +  p ( 7, 12) $   and
$\bar{\t} \cap \Z^2 = \bigsqcup_{(p, q) \in \Z^2_{\geq 0}}  G' + p (1,1) +  q (7,12)$.
\end{center}
We deduce that
$F_{\bar{A}}   = F_{\bar{\t} \cap \Z^2}  - F_{ (\bar{\t} \cap \Z^2 )\setminus \bar{A}}$ hence
\begin{center}
$F_{\bar{A}} = (1- x_1^7 x_2 ^{12})^{-1}  \left( \sum_{(i,j) \in G} x_1^i x_2^j +
(\sum_{(i,j) \in G'} x_1^i x_2^j) (1 - x_1 x_2 )^{-1}   \right) $.
\end{center}
               To get the series $P_{1,\t}(\Lambda)$ we apply to
$F_{\bar{A}} $ the ring homomorphism which maps
$x_1^i x_2^j \mapsto \L^{j-i} T^j$ and then  we multiply the result by
$\L-1$.

 We check that none of the candidate
poles of $P_{\mathrm {geom}}^{(Z^\Lambda,0)}(T)$ cancels. The
motivic volume  is
\begin{center}
$\mu(H_\Lambda)=(\L-1)^2\Big(\frac{1}{(1-\L)(1-\L^{19})}+
\frac{-1}{1-\L^{19}}+\frac{1+\L^8+\L^{16}}{(1-\L^{19})(1-\L^5)}+\frac{-1}{1-\L^5}+\frac{1}{(1-\L^5)(1-\L)}\Big)$.
\end{center}

         \vspace{0.3cm}

\noindent {\bf Acknowledgments:} We are grateful to Monique
Lejeune-Jalabert and to Johannes Nicaise for their suggestions,
comments and motivations for this work, and to the Section of
Algebra of the Mathematics Department of the K.U. Leuven for the
hospitality.
 
\end{document}